\documentclass[10pt]{article}

\usepackage{amsfonts}
\usepackage{amsmath}
\usepackage{amssymb}
\usepackage{enumerate}
\usepackage{framed}
\usepackage{graphicx}
\usepackage{lscape}
\usepackage{bm}
\usepackage{color}
\usepackage{multicol}

\usepackage{amsthm} 
\newtheorem{theorem}{Theorem}[section]
\newtheorem{corollary}[theorem]{Corollary}
\newtheorem{lemma}[theorem]{Lemma}
\newtheorem{proposition}[theorem]{Proposition}
\newtheorem{definition}[theorem]{Definition}
\newtheorem{remark}[theorem]{Remark}

\newenvironment{pro}{\vspace{1ex}\noindent{\bf Proof:}\hspace{0.5em}}
	{\hfill\qed\vspace{1ex}}

\newcommand{\R}{\mathbb{R}}
\newcommand{\Z}{\mathbb{Z}}
\newcommand{\N}{\mathbb{N}}

\newcommand{\C}{\mathbb{C}}
\newcommand{\D}{\partial}

\newcommand{\w}{\omega}
\newcommand{\Y}{\lambda}

\newcommand{\s}{\sigma}
\newcommand{\op}{\oplus}
\newcommand{\bu}{\bm{u}}
\newcommand{\bv}{\bm{v}}

\title{\textbf{Abelian functions associated with a cyclic tetragonal curve of genus six}}
\author{M. ENGLAND AND J.C. EILBECK \\ Heriot Watt University}
\date{\today}

\begin{document}

\maketitle

\subsubsection*{Abstract}
We develop the theory of Abelian functions defined using a tetragonal curve of genus six, discussing in detail the cyclic curve $y^4 = x^5 + \lambda_4x^4 + \lambda_3x^3 + \lambda_2x^2 + \lambda_1x + \lambda_0$.  We construct Abelian functions using the multivariate $\sigma$-function associated to the curve, generalising the theory of the Weierstrass $\wp$-function.  We demonstrate that such functions can give a solution to the KP-equation, outlining how a general class of solutions could be generated using a wider class of curves.  We also present the associated partial differential equations satisfied by the functions, the solution of the Jacobi Inversion Problem, a power series expansion for $\sigma(\bu)$ and a new addition formula.

\section{Introduction} \label{SECintro}

Recent times have seen a revival of interest in the theory of Abelian (multiply periodic) functions associated with algebraic curves.  This topic can be dated back to the Weierstrass theory of elliptic functions, which we use as a model.  Let $\sigma(u)$ and $\wp(u)$ be the standard Weierstrass functions (see for example \cite{WW47}).  The $\wp$-function can be used to parametrise an elliptic curve, $y^2 = 4x^3 - g_2x - g_3$, and satisfies the following well-known formulae.
\begin{align}
      \wp(u) &= - \frac{d^2}{du^2} \log \sigma(u),      \label{elliptic_ps}     \\
\big(\wp'(u)\big)^2 &= 4\wp(u)^3 - g_2\wp(u) - g_3,  \label{elliptic_diff1}     \\
    \wp''(u) &= 6\wp(u)^2 - \tfrac{1}{2}g_2.         \label{elliptic_diff2}
\end{align}
The $\sigma$-function satisfied a power series expansion
\begin{align}
\sigma(u) &= u - \frac{1}{240}g_2u^5 - \frac{1}{840}g_3u^7 - \frac{1}{161280}g_2^2u^9 
- \frac{1}{2217600}g_2g_3u^{11} + \dots                                       \label{elliptic_sigexp}\
\end{align}
and a two term addition formula 
\begin{align}
- \frac{\sigma(u+v)\sigma(u-v)}{\sigma(u)^2\sigma(v)^2} = \wp(u) - \wp(v).  \label{elliptic_add}
\end{align}
Taking logarithmic derivatives of this will give the standard addition formula for $\wp(u)$. This paper will generalise equations (\ref{elliptic_ps})-(\ref{elliptic_add}) for a previously unconsidered class of functions..

\vspace{0.05in}

The study of Abelian functions associated with the simplest hyperelliptic curves (those of genus 2) goes back to the start of the 20th century.  Klein's generalisation of the Weierstrass theory is described in Baker's classic texts, \cite{ba97} and \cite{ba07}, while \cite{bel97} gives a more recent study of the general hyperelliptic case.
Further generalisation has been structured by considering, (with the notation of \cite{bc23}), classes of \emph{$(n,s)$-curves}.  These are curves with equation
\begin{equation} \label{ns_curve}
y^n - x^s - \sum_{\alpha,\beta} \mu_{[ns - \alpha n - \beta s]} x^{\alpha}y^{\beta} \qquad \mu_j \mbox{ constants,}
\end{equation}
where $\alpha, \beta \in \Z$ with $\alpha \in (0,s-1), \beta \in (0,n-1)$ and $\alpha n + \beta s < ns$.  The cyclic subset of such a class of curves is generated by setting $\beta=0$.  We suppose that $(n,s)$ are coprime, in which case the curves have genus $g = \frac{1}{2}(n-1)(s-1)$, and a unique branch point $\infty$ at infinity.

In the last few years a good deal of progress has been made on the theory of Abelian functions associated with trigonal curves, (those with $n=3$). The $\sigma$-function realisation of these functions was developed first in \cite{bel00} and \cite{eel00}, with the two canonical cases studied in detail in \cite{eemop07} and \cite{bego08}. 

\vspace{0.05in}

In this paper we consider the next logical class, and work with a tetragonal curve, (one with $n=4$).  We have started by looking at the curves of lowest genus, and simplified by considering the cyclic subclass.  We construct the multivariate $\sigma$-function associated with this curve, and use it to define and analyse classes of Abelian functions, generalising the theory of the Weierstrass $\wp$-function.  A key component of our work was the construction of a series expansion for the $\sigma$-function.  This technique was first developed for the trigonal case in \cite{bg06}, however the computation involved for the present expansion is significantly greater.  The latter computations were performed in parallel with the use of the Distributed Maple software, (see \cite{DM} and \cite{smb03}).

The applications of Abelian functions to integrable systems and soliton theory have been the topic of research for some time, (see for example \cite{DoublePend} and \cite{EEP01}).  It is well known that the elliptic $\wp$-function could be used to construct a solution to the KdV-equation.  Similar solutions to non-linear equations have been derived from higher genus curves, for example in \cite{eel00}, where the function $\wp_{33}$ associated with the (3,4)-curve was shown to be a solution of the Boussinesq equation.  This has suggested a more general link between such functions and the integrable KP hierarchy.  We have demonstrated how the Abelian functions we define can give a solution to the KP-equation, outlining how similar solutions will also be found from any $(4,s)$-curve.  

This paper is organised as follows.  We give the basic properties of the curve we consider in Section \ref{SECcurvedef}, including explicit constructions of the differentials on the curve, and a set of weights that render the key equations homogeneous.  Then in Sections \ref{SECsigdef} and \ref{SECabeliandef} we define the $\sigma$-function and Abelian functions associated with this curve.  Section \ref{SECkleinform} discusses a key theorem satisfied by the $\wp$-functions, which we use to give a solution to the Jacobi Inversion Problem.  
In Section \ref{SECsigexp} we derive some properties of the $\sigma$-function, including the series expansion, while in Section \ref{SECbasis} we use this to generate relations between the Abelian functions.  Section \ref{SECnonlin} demonstrates how solutions to the KP-equation can be constructed from Abelian functions.  Finally, in Section \ref{SECaddform} we give the derivation of a two-term addition formula.

\section{The purely tetragonal curves} \label{SECcurvedef}

We will investigate Abelian functions associated with a tetragonal curve.  The simplest general tetragonal curve is, in the notation of the $(n,s)$-curves, a (4,5)-curve.  This would be given by $g(x,y)=0$ where
\begin{align*}
&g(x,y) = y^4 + \big( \mu_{1}x + \mu_{5} \big)y^3
+ \big( \mu_{2}x^2 + \mu_{6}x  + \mu_{10} \big)y^2 \nonumber \\
&\quad + \big( \mu_{3}x^3 + \mu_{7}x^2 + \mu_{11}x + \mu_{15} \big)y \nonumber \\
&\quad - \big( x^5 + \mu_{4}x^4 + \mu_{8}x^3 + \mu_{12}x^2 + \mu_{16}x + \mu_{20} \big) \qquad 
(\mu_j \mbox{ constants).}
\end{align*}
In this paper we further simplify by considering the \emph{cyclic subclass} of this family.  These are the curves $C$, given by 
\begin{align}
C: \quad &f(x,y) = 0  
\hspace*{2in} (\lambda_j \mbox{ constants)} \nonumber \\
\qquad \mbox{where} \quad &f(x,y) = y^4 - \big( x^5 + \lambda_4 x^4 + \lambda_3 x^3 + \lambda_2 x^2+\lambda_1 x + \lambda_0 \big).\label{C} 
\end{align}
The curve $C$ has genus $g=6$, the unique branch point $\infty$ at infinity and is referred to as the \textit{purely tetragonal}, or \textit{strictly tetragonal} curve. It contains an extra level of symmetry, demonstrated by the fact that it is invariant under
\begin{equation} \label{invar}
[\zeta]: \quad (x,y) \rightarrow (x, \zeta y),
\end{equation}
where $\zeta$ is a $4$th root of unity.

For any $(n,s)$-curve we can define a set of weights for the variables of the theory, including the curve constants, which render equation homogeneous with respect to the weights.  To find these weights consider the mapping 
$\chi \mapsto t^{\alpha_{\chi}} \tilde{\chi}$ acting on all elements in the curve equation.  Then define the weights as the constants $\alpha_i$ that render the new equation homogeneous with respect to $t$.  The weights of $x,y$ will then be determined up to a constant by $n\alpha_y = s\alpha_x$.  To keep with convention, we let $\alpha_{x}=-n$ and $\alpha_y=-s$ so that they are the largest negative integers satisfying this condition.  The weights of the curve constants can then be chosen to make the remainder of the equation homogeneous.  

\begin{definition}\label{Sato1}
For the cyclic (4,5)-case we have 
\begin{center}
\begin{tabular}{|c|c|c|c|c|c|c|c|c|c|c|c|c|}\hline
                & $x$  & $y$   &  $\lambda_4$ & $\lambda_3$ & $\lambda_2$ & $\lambda_1$ & $\lambda_0$    \\ \hline
\textbf{Weight} & $-4$ & $-5$  &  $-4$        & $-8$        & $-12$       & $-16$       & $-20$          \\ \hline
\end{tabular}
\end{center}
while in the general (4,5)-case the weights of the curve constants are given by their subscripts.  We refer to these as the \emph{Sato Weights}
\end{definition}
\noindent As we procede through the paper we can use the approach of this mapping to conclude that other elements in our theory must have definite weight, and all the equations presented here will be homogeneous with respect to these weights. 

Next we construct the standard basis of holomorphic differentials upon $C$. 
\begin{align}
\bm{du} &= (du_1, \dots ,du_6), \qquad du_i(x,y) = \frac{g_i(x,y)}{4y^3}dx, \nonumber \\
&\mbox{where} \qquad 
\begin{array}{lllll}
g_1(x,y) = 1,   & \quad & g_2(x,y) = x, & \quad & g_3(x,y) = y, \\
g_4(x,y) = x^2, & \quad & g_5(x,y) =xy, & \quad & g_6(x,y) = y^2.\label{holodiff}
\end{array}
\end{align}
Denote points in $\C^6$ by $\bu$ for example, and their coordinates by $(u_1,u_2, \dots ,u_6)$.  We know from the general theory, that any point $\bu \in \C^6$ can be expressed as 
\begin{align*}
\bu &= (u_1,u_2,u_3,u_4,u_5,u_6) = \sum_{i=1}^6 \int_{\infty}^{P_i} \bm{du},
\end{align*}
where the $P_i$ are six variable points upon $C$. Let $\Lambda$ denote the lattice generated by the integrals of the basis of holomorphic differentials along any closed paths in $C$.  Then the manifold $\C^6 / \Lambda$ is the Jacobian variety of $C$, denoted by $J$.  Let $\kappa$ be the map of modulo $\Lambda$ over $\C$:
\begin{equation*}
\kappa: \C^6 \to \C^6 / \Lambda = J.
\end{equation*}
Therefore $\Lambda = \kappa^{-1}((0,\dots,0))$.  Next, for $k=1,2,\dots$ define $\mathfrak{A}$, the \emph{Abel map} from the $k$th symmetric product Sym$^k(C)$ to $J$.
\begin{align}
\mathfrak{A}: \mbox{Sym}^k(C) &\to     J \nonumber \\
(P_1,\dots,P_k)   &\mapsto \left( \int_{\infty}^{P_1} \bm{du} + \dots + \int_{\infty}^{P_k} \bm{du} \right) \pmod{\Lambda}.
\label{Abel}
\end{align}
where the $P_i$ are points upon $C$.  Denote the image of the $k$th Abel map by $W^{[k]}$, and let
\begin{equation*}
[-1](u_1, \dots ,u_6) = (-u_1, \dots ,-u_6).
\end{equation*}
Define the \emph{$k$th standard theta subset} (often referred to as the $k$th strata) by
\begin{equation*}
\Theta^{[k]} = W^{[k]} \cup [-1]W^{[k]}.
\end{equation*}
When $k=1$ the Abel map gives an embedding of the curve $C$, upon which we define $\xi$ as the local parameter at the origin, $\mathfrak{A}_1(\infty)$.
\begin{equation} \label{xi_def}
\xi = x^{-\frac{1}{n}} = x^{-\frac{1}{4}}
\end{equation}
We can then express the basis (\ref{holodiff}) with $\xi$ and integrate to give
\begin{align}
\begin{array}{lll}
u_1 = -\frac{1}{11}\xi^{11} + O(\xi^{15}) & u_3= -\frac{1}{6}\xi^6 + O(\xi^{10}) & u_5 = -\frac{1}{2}\xi^2 + O(\xi^6) \\
u_2 = -\frac{1}{7}\xi^{7} + O(\xi^{11})   & u_4= -\frac{1}{3}\xi^3 + O(\xi^7)    & u_6 = -\xi + O(\xi^5).
\end{array} \label{diffs_xi}
\end{align}
The higher order terms will contain the curve parameters $\bm{\lambda} = \{\lambda_0, \dots, \lambda_4\}$.
Similarly to \cite{bg06} and \cite{eel00}, we could rewrite this using $u_6$ as the local paremeter.  

Note that such calculations can be performed similarly for any $(n,s)$-curve, and that since each element of $\bm{du}$ is homogeneous in Sato weight we can conclude that the $u_i$ have definite Sato weight.  Since the weight of $\xi$ must be $+1$ from equation (\ref{xi_def}) we can define the weights of $\bm{u}$ uniquely as below.
\begin{definition}\label{Sato2}
In the (4,5)-case we assign the following weights to $\bu$. 
\begin{center}
\begin{tabular}{|c|c|c|c|c|c|c|c|c|c|c|c|}\hline
                & $u_1$ & $u_2$  &  $u_3$ & $u_4$ & $u_5$ & $u_1$ \\ \hline
\textbf{Weight} & $+11$ & $+7$   &  $+6$  & $+3$  & $+2$  & $+1$  \\ \hline
\end{tabular}
\end{center}
\end{definition}
\begin{remark}
The weights of the variables coincide with the order of their zero at $\infty$.  They can also be calculated using the Weierstrass gap sequence, where the weights of $u_1,\dots,u_6$ are the gap numbers, and the weights of $x$ and $y$ are the negative of the first two non-gap numbers.
\end{remark}

\begin{definition}
Let $(x,y)$ and $(z,w)$ be two variable points upon $C$.  Then the 2-form $\Omega \big( (x,y), (z,w) \big)$ on $C \times C$ is a \emph{fundamental differential of the second kind} if 
\begin{enumerate}[1.]
\item It is symmetric: $\quad \Omega \big( (x,y), (z,w) \big) = \Omega \big( (z,w), (x,y) \big)$.
\item The only pole of second order is along the diagonal of $C \times C$ (where $x=z$).
\item It can be expanded in a power series as
\[
\Omega \big( (x,y), (z,w) \big) = \left( \frac{1}{(\xi - \xi')^2} + O(1) \right) d \xi d \xi' \quad (\mbox{as } (x,y) \to (z,w)),
\]
where $\xi$ and $\xi'$ are the local coordinates of $(x,y)$ and $(z,w)$.
\end{enumerate}
\end{definition}
\noindent We will construct Klein's explicit realisation of this in Proposition \ref{Omegexp} below.

First introduce $\bm{dr}$, the basis of meromorphic differentials which have their only pole at $\infty$. These are determined modulo the space spanned by the $\bm{du}$ and can be expressed as
\begin{equation} \label{merodiffbas}
\bm{dr} = (dr_1,\dots,dr_6), \qquad \mbox{where} \quad dr_j(x,y) = \frac{h_j(x,y)}{4y^3}dx.
\end{equation}
An explicit basis is constructed later, in order to satisfy Proposition \ref{Omegexp}.
\begin{definition} Define the following meromorphic function on $C \times C$ as
\begin{equation*}
\Sigma \big( (x,y), (z,w) \big) = \frac{1}{4y^3(x-z)} \cdot \sum_{k=1}^4 y^{4-k} \left[ \frac{f(z,w)}{w^{4-k+1}} \right]_w ,
\end{equation*}
where $[\quad]_w$ means that we remove any terms which have negative powers with respect to $w$.
\end{definition}
\begin{proposition} \label{Omegexp}
The fundamental differential of the second kind can be expressed as 
\begin{equation*}
\Omega \big( (x,y), (z,w) \big) = R \big( (x,y), (z,w) \big)dxdz,
\end{equation*}
where 
\begin{equation*}
R\big( (x,y), (z,w) \big) = \frac{\partial}{\partial z} \Sigma \big( (x,y), (z,w) \big) 
+ \sum_{j=1}^6 \frac{du_j(x,y)}{dx} \cdot \frac{dr_j(z,w)}{dz}.
\end{equation*}
The polynomials $h_j(x,y)$ need to be chosen so that $\Omega$ is symmetric.  This will lead to a realisation of $\Omega$ in the form
\begin{equation} \label{OmgF}
\Omega \big( (x,y), (z,w) \big) = \frac{ F \big( (x,y), (z,w) \big) dx dz}{(x-z)^2f_y(x,y)f_w(x,y)}.
\end{equation}
\end{proposition}
\begin{pro}
The essential part of the proof is the same as in the lower genus cases (see \cite{eemop07} for example).  In this case we explicitly determine the basis of meromorphic differentials, (\ref{merodiffbas}), to be given with 
\begin{align*}
h_1 &= -y^2 \big( 8x^2\lambda_4 + 11x^3 + 5x\lambda_3 + 2\lambda_2 \big), \qquad
h_2 = -y^2 \big( \lambda_3 + 4x\lambda_4 + 7x^2 \big), \nonumber \\
h_3 &= -2xy \big( \lambda_3 + 3x^2 + 2x\lambda_4 \big), \quad
h_4 = -3xy^2, \quad
h_5 = -2x^2y, \quad
h_6 = -x^3.
\end{align*}
The polynomial $F$ in the realisation (\ref{OmgF}) is found to be,
\begin{align}
&F \big( (x,y), (z,w) \big) = 
4y^3w^3 + \big( 3xz^4 + z^3\lambda_3 + z^3x^2 + 2\lambda_2z^2 + 3x\lambda_3z^2 \nonumber \\
&\quad + 4z^3x\lambda_4 + 4\lambda_0 + \lambda_1x + 2\lambda_2xz + 3\lambda_1z \big)y^2 
+ \big( 2\lambda_1z + 4\lambda_2xz + 4\lambda_0 \nonumber \\
&\quad + 2\lambda_1x + 4x^2\lambda_4z^2 + 2\lambda_3x^2z + 2x^3z^2 
+ 2z^3x^2 + 2x\lambda_3z^2 \big) wy \nonumber \\
&\quad + \big( \lambda_3x^3 + 4\lambda_0 + 3\lambda_1x + 2\lambda_2x^2 + \lambda_1z + x^3z^2 + 3x^4z + 2\lambda_2xz 
\nonumber \\
&\quad + 3\lambda_3x^2z + 4\lambda_4x^3z \big) w^2. \label{fundF}
\end{align}
\end{pro}

\section{Defining the $\bm{\sigma}$-function}\label{SECsigdef}

In this section we describe the multivariate $\sigma$-function associated with $C$, from which all Abelian functions associated with $C$ can be defined.  This can be regarded as a generalisation of the Weierstrass $\sigma$-function, with the main difference that there are now $g=6$ variables.
\begin{equation*}
\sigma = \sigma(\bu) = \sigma(u_1,u_2,u_3,u_4,u_5,u_6).
\end{equation*}
First we choose a basis of cycles (closed paths) upon the surface defined by $C$.  We denote them
\begin{equation*}
\alpha_i, \beta_j, \qquad 1 \leq i,j \leq 6,
\end{equation*}
and ensure they have intersection numbers
\begin{equation*}
\alpha_i \cdot \alpha_j = 0, \qquad \beta_i \cdot \beta_j = 0, \qquad 
\alpha_i \cdot \beta_j = \delta_{ij} = \left\{  
\begin{array}{ccc}
1 & \mbox{if} & i = j \\
0 & \mbox{if} & i \neq j
\end{array}.  \color{white} \right\} 
\end{equation*}
This allows us to define the following period matrices:
\begin{align*}
\begin{array}{cc}
          \w'  = \left( \oint_{\alpha_k} du_\ell \right)_{k,\ell = 1,\dots,6} &
\qquad    \w'' = \left( \oint_{ \beta_k} du_\ell \right)_{k,\ell = 1,\dots,6}  \\
        \eta'  = \left( \oint_{\alpha_k} dr_{\ell} \right)_{k,\ell = 1,\dots,6} &
\qquad  \eta'' = \left( \oint_{ \beta_k} dr_{\ell} \right)_{k,\ell = 1,\dots,6}
\end{array}.
\end{align*}
We combine these into
\begin{equation*}
M =
\begin{bmatrix}
\w' & \w'' \\
\eta' & \eta''
\end{bmatrix},
\end{equation*}
which we know from classical results to satisfy
\begin{equation}\label{Legendre}
M 
\begin{bmatrix}
  & -I_6 \\
I_6  
\end{bmatrix}
^T M = 2\pi i
\begin{bmatrix}
\w' & \w'' \\
\eta' & \eta''
\end{bmatrix}.
\end{equation}
This is the \emph{generalised Legendre equation}, (see \cite{bel97} p11).  We also have that $(\omega')^{-1}\omega''$ is symmetric with 
\begin{equation} \label{imomeg}
\mbox{Im} \big( (\omega')^{-1}\omega'' \big) \quad \mbox{ positive definite}.
\end{equation}
We now define the multivariate $\sigma$-function associated with $C$.  This can be constructed using the multivariate $\theta$-function, (see for example, \cite{mu83}).
\begin{definition} \label{sigdef}
The \emph{Kleinian $\s$-function associated with $C$} is
\begin{align*}
\s(\textbf{u}) &\underset{\textcolor{white}{y}}{=} \sigma(\bu; M) = c \exp \big( - \textstyle \frac{1}{2} \bm{u} \eta' (\w')^{-1} \bm{u}^T \big) \times \theta[\delta]\big((\w')^{-1}\bm{u}^T \hspace*{0.05in} \big| \hspace*{0.05in} (\w')^{-1} \w''\big) \nonumber \\
 &= c \exp \big( - \textstyle \frac{1}{2} \bm{u} \eta' (\w')^{-1} \bm{u}^T \big) \\
\times \sum_{m \in \Z^6} &\exp \bigg[ 2\pi i \bigg\{ \textstyle \frac{1}{2} (m+\delta')^T (\w')^{-1} \w''(m+\delta') + 
(m+\delta')^T ((\w')^{-1} \bm{u}^T + \delta'') \bigg\} \bigg], \nonumber
\end{align*}
where $c$ is a constant dependent upon the curve parameters, $\{\lambda_0, \lambda_1, \lambda_2, \lambda_3, \lambda_4\}$ and fixed later (see Remark \ref{remfixc}).  The matrix 
$\delta = \left[ \begin{array}{l}
\bm{\delta'} \\ \bm{\delta''}
\end{array}
\right]$
is the theta function characteristic which gives the Riemann constant for $C$ with respect to the base point $\infty$ and the period matrix $[\w', \w'']$, (see \cite{bel97} p23-24).
\end{definition}
In this paper we give some of the most important properties of $\sigma(\bu)$.  However, for a more detailed study of the construction and properties of the multivariate $\sigma$-function, we refer the reader to \cite{bel97}.
\begin{lemma} \label{sigper}
Given $\bu \in \C^6$, denote by $\bm{u'}$ and $\bm{u''}$ the unique elements in $\R^6$ such that 
\begin{equation*}
\bu = \bm{u'}\w' + \bm{u''}\w''.
\end{equation*}
Let $\ell$ represent a point on the period lattice
\begin{equation*}
\ell = \ell'\w' + \ell''\w'' \in \Lambda.
\end{equation*}
For $\bu, \bv \in \C^6$ and $\ell \in \Lambda$, define $L(\bu,\bv)$ and $\chi(\ell)$ as follows:
\begin{align*}
L(\bu,\bv) &= \bu^T \big( \eta'\bm{v'} + \eta''\bm{v''} \big), \\
\chi(\ell) &= \exp \big[ \pi i \big( 2(\ell'^T\delta'' - \ell''^T\delta') + \ell'^T\ell'' \big) \big].
\end{align*}
Then, for all $\bu \in \C^6, \ell \in \Lambda$ the function $\sigma(\bu)$ has the quasi-periodicity property
\begin{equation} \label{quas}
\sigma(\bu + \ell) = \chi(\ell) \exp \Big[ L \Big( u + \frac{\ell}{2}, \ell \Big) \Big] \cdot \sigma(\bu).
\end{equation}
Also, for $\gamma \in Sp(12,\Z)$ we have
\begin{equation} \label{Mper}
\sigma(\bu; \gamma M) = \sigma(\bu; M).
\end{equation}
\end{lemma}
\begin{pro}
The quasi-periodicity property given in equation (\ref{quas}) is a classical result, first discussed in \cite{ba97}, that was fundamental to the original definition of the multivariate $\sigma$-function.  Equation (\ref{Mper}) is easily seen from the definition of $\sigma(\bu)$, since $\gamma$ corresponds to the choice of basis cycles $\{\alpha_j, \beta_j \}_{j=1}^6$ which were used to define $M$. \\
\end{pro}

\section{Classes of Abelian functions}\label{SECabeliandef}

\begin{definition}
Let $\mathfrak{M}(\bu)$ be a meromorphic function of $\bu \in \C^6$.  Then $\mathfrak{M}$ is an \emph{Abelian function associated with $C$} if
\begin{equation*}
\mathfrak{M}(\bu + \w' \bm{n}^T + \w'' \bm{m}^T) = \mathfrak{M}(\bu),
\end{equation*}
for all integer vectors $\bm{n},\bm{m} \in \Z$,  wherever $\mathfrak{M}$ is defined.
\end{definition}
\noindent We now define a set of fundamental Abelian functions on $J$.
\begin{definition} \label{defnip}
Define the \emph{2-index Kleinian $\wp$-functions} as
\begin{equation*}
\wp_{ij}(\bu) = - \frac{\D^2}{\D u_i\D u_j} \log \sigma(\bu), \qquad i \leq j \in \{1,\dots,6\}.
\end{equation*}
\end{definition}
\noindent A short calculation shows these functions to have have poles of order 2 when $\sigma(\bu)=0$, and no other singularities.  We can check (using Lemma \ref{sigper}), that
\[
\wp_{ij}(\bu + \ell) = \wp_{ij}(\bu), \qquad \forall \quad \ell \in \Lambda.
\] 
Hence we can conclude these functions to be Abelian.  Similar analysis will show their derivatives to be Abelian also.
\begin{definition} \label{nordp} For $n \geq 2$, define \emph{$n$-index Kleinian $\wp$-functions} as
\begin{equation*}
\wp_{i_1,i_2,\dots,i_n}(\bu) = - \frac{\D}{\D u_{i_1}} \frac{\D}{\D u_{i_2}}\dots \frac{\D}{\D u_{i_n}} \log \sigma(\bu),
\quad i_1 \leq \dots \leq i_n \in \{1,\dots,6\}. 
\end{equation*}
\end{definition}
\begin{remark}
(i) Compare with equation (\ref{elliptic_ps}) and the elliptic case to see that we are defining a generalisation of the Weierstrass $\wp$-function, and its derivatives. \\
(ii) This notation is compatible with the elliptic case, where we would now denote the Weierstrass $\wp$-function as $\wp_{11}(\bu)$ and its first derivative $\wp'$ by $\wp_{111}(\bu)$.  \\
(iii) The order of the indices is irrelevant.  For simplicity we always use ascending numerical order. \\
(iv) We are usually only referring to one vector of variables $\bu$.  
In these cases, for simplicity, we write $\wp_{ij}$ instead of $\wp_{ij}(\bu)$.
\end{remark}

\noindent We find in Section \ref{SECbasis} that the $\wp$-functions are not sufficient to construct a basis of the simplest Abelian functions.  Hence we also define a generalisation of Baker's $Q$-functions, which we need to extend further than in the lower genus cases. 
\begin{definition} \label{Qdef}
Define the operator $\Delta_i$ as below.  This is now known as \emph{Hirota's bilinear operator}, although it was used much earlier by Baker in \cite{ba07} for example.
\[
\Delta_i = \dfrac{\D}{\D u_i} - \dfrac{\D}{\D v_i}. 
\]
It is then simple to check that an alternative, equivalent definition of the 2-index Kleinian $\wp$-functions is given by
\begin{align*}
\wp_{ij}(\bm{u}) &= - \frac{1}{2\s(\bm{u})^2} \Delta_i\Delta_j \s(\bm{u}) \s(\bm{v}) \hspace*{0.1in} \Big|_{v=u}
\qquad i \leq j \in \{1,\dots,6\}.
\end{align*}
We extend this to define \emph{$n$-index $Q$-functions}, for $n$ even.
\[
Q_{i_1, i_2,\dots,i_n}(\bm{u}) =  \frac{(-1)}{2\s(\bm{u})^2} \Delta_{i_1}\Delta_{i_2}...\Delta_{i_n} \s(\bm{u}) \s(\bm{v}) \hspace*{0.08in} \Big|_{v=u} \quad i_1 \leq ... \leq i_n \in \{1,\dots,6\}.
\]
\end{definition}
\noindent We can show, as above, that these functions are also Abelian.

\begin{remark}\label{Qspec}
(i) The subscripts of the $\wp$-functions denote differentiation
\[
\frac{\D}{\D u_{i_n+1}}\wp_{i_1,i_2,\dots,i_n} = \wp_{i_1,i_2,\dots,i_n, i_{n+1}},
\]
but this is not the case for the $Q$-functions.  Here the indices refer to which Hirota operators were used. \\
(ii) If we had applied the definition for $n$ odd, then it would have returned zero. \\
(iii) Note that both the $\wp$-functions and the $Q$-function have poles when $\sigma(\bu)=0$ and no other singularities.  However, the $n$-index $\wp$-functions had poles of order $n$, while the $n$-index $Q$-function all have poles of order 2.
\end{remark}
The 4-index $Q$-functions were first used by Baker, and in \cite{eemop07} it was shown that they could be expressed using the Kleinian $\wp$-functions as 
\begin{equation} \label{4iQ}
Q_{ijk\ell} = \wp_{ijk\ell} - 2 \wp_{ij}\wp_{k\ell}-2\wp_{ik}\wp_{j\ell} -2\wp_{i\ell}\wp_{jk}.
\end{equation}
\begin{proposition} 
The 6-index $Q$-functions can be written as
\begin{align}
Q_{ijklmn} &= \wp_{ijklmn} - 2\Big( 
\big[\wp_{ij}\wp_{klmn} + \wp_{ik}\wp_{jlmn} + \wp_{il}\wp_{jkmn} + \wp_{im}\wp_{jkln} \nonumber \\
&\quad + \wp_{in}\wp_{jklm}\big] + \big[ \wp_{jk}\wp_{ilmn} + \wp_{jl}\wp_{ikmn} + \wp_{jm}\wp_{ikln} + \wp_{jn}\wp_{iklm} \big] \nonumber \\
&\quad + \big[ \wp_{kl}\wp_{ijmn} + \wp_{km}\wp_{ijln} + \wp_{kn}\wp_{ijlm} \big] + \big[ \wp_{lm}\wp_{ijkn} + \wp_{ln}\wp_{ijkm} \big] \nonumber \\
&\quad + \wp_{mn}\wp_{ijkl} \Big) + 4\Big( \big[ \wp_{ij}\wp_{kl}\wp_{mn} + \wp_{ij}\wp_{km}\wp_{ln} + \wp_{ij}\wp_{kn}\wp_{lm} \big] \nonumber \\
&\quad + \big[\wp_{ik}\wp_{jl}\wp_{mn} + \wp_{ik}\wp_{jm}\wp_{ln} + \wp_{ik}\wp_{jn}\wp_{lm} \big] 
+ \big[\wp_{il}\wp_{jk}\wp_{mn} \nonumber \\
&\quad + \wp_{il}\wp_{jm}\wp_{kn} + \wp_{il}\wp_{jn}\wp_{km} \big] + \big[\wp_{im}\wp_{jk}\wp_{ln} 
+ \wp_{im}\wp_{jl}\wp_{kn} \nonumber \\
&\quad + \wp_{im}\wp_{jn}\wp_{kl} \big] + \big[\wp_{in}\wp_{jk}\wp_{lm} + \wp_{in}\wp_{jl}\wp_{km} 
+ \wp_{in}\wp_{jm}\wp_{kl} \big] \Big). \label{6iQ}
\end{align}
\end{proposition}
\begin{pro}
Apply Definitions \ref{defnip} and {Qdef} to reduce the equation to a sum of $\sigma$-derivatives.  We find that they all cancel (Maple is useful here).  The structure of the sum was prompted by considering the result for the 4-index $Q$-functions.
\end{pro}

\noindent Clearly, equation (\ref{6iQ}) will specialise to give a set of simpler formulae, such as 
\begin{align*}
Q_{nnnnnn} = \wp_{nnnnnn} - 30\wp_{nn}\wp_{nnnn} + 60\wp_{nn}^3.
\end{align*}

\section{Expanding the Kleinian formula}\label{SECkleinform}

This section is based upon the following Theorem (originally by Klein).  It is given for a general curve as Theorem 3.4 in \cite{eel00}.  From this theorem we are able to solve the Jacobi Inversion Problem, as well as generate relations between the $\wp$-functions 
\begin{theorem} \label{Klein_Thm}
Let $\{P_1,...,P_6\} \in C^6$ be an arbitrary set of distinct points on $C$, and $(z,w)$ any point of this set.  Then for an arbitrary point $(x,y)$ and base point $\infty$ on $C$ we have
\begin{equation} \label{Kleineq}
\sum_{i,j=1}^6 \wp_{ij} \left( \int_{\infty}^t d\bm{u} - \sum_{k=1}^6 \int_{\infty}^{P_k} d\bm{u} \right)
g_i(x,y) g_j(z,w)
= \frac{F\big((x,y),(z,w)\big)}{(x-z)^2}.
\end{equation}
Here $g_i$ is the numerator of $du_i$, as given in equation (\ref{holodiff}), and $F$ is the symmetric function appearing in equation (\ref{fundF}) as the numerator of the fundamental differential of the second kind.
\end{theorem}
We use our explicit calculation of the differentials to construct (\ref{Kleineq}).  We expand this as one of the $P_k$ tends to infinity, to obtain a series expansion in terms of the local parameter $\xi$, given earlier in equation (\ref{xi_def}).  It follows that each coefficient with respect to $\xi$ must be zero for any $\bm{u} \in J$ and some $(z,w)$ on $C$.  
This gives us a potentially infinite sequence of equations, starting with the five given in Appendix \ref{APP_KE}.  The first 14 have been calculated explicitly (using Maple), and can be found online at \cite{myweb}.

\subsubsection*{Manipulating these equations}

We follow the approach of the trigonal papers, and manipulate these equations using Maple.  We first take the resultant of pairs of these equations, (eliminating the variable $w$ by choice), to give a new set of equations dependent on $z$ and the $\wp$-functions.
We let  Res($a,b$) represent the resultant of equations ($a$) and ($b$). 

These new equations are considerably longer than those obtained in the lower genus cases. 
We need to combine them to give a polynomial of degree $g-1=5$ in $z$.  Such a polynomial would have only 5 solutions, but must be satisfied for all $\bu$ (which has 6 variables).  Hence all the coefficients must be zero, giving us a set of relations between the $\wp$-functions.  However, we have the extra complication, (compared to the trigonal cases), that none of the new equations has degree in $z$ equal to $g$.  Therefore, in this case, at least two rounds of elimination between the equations will be required.  

We find that Res($\ref{pp1},\ref{pp2}$) has degree 7 in z, so we rearrange it to give an equation for $z^7$.  Then since Res($\ref{pp1},\ref{pp3}$) and Res($\ref{pp1},\ref{pp5}$) have degree 8, we can repeatedly substitute for $z^7$ in both until we are left with two equations of degree 6 in $z$.  Since these are very long we do not print them here, however they can be found online at \cite{myweb} where we have labelled them (T1) and (T2).

We next rearrange (T1) to give an equation for $z^6$ and repeatedly substitute for $z^7$ and $z^6$ in the remaining equations until they are of degree 5 in $z$.  The coefficients of such equations must be zero, giving us relations between the $\wp$-functions.   
The smallest such relation has 3695 terms, with the others rising in size considerably.  Unlike the trigonal cases, these can not be easily separated to give expressions for individual $\wp$-functions.  However, these are implemented in the construction of the $\sigma$-function expansion, (see Section \ref{SECsigexp}).

\subsubsection*{Jacobi Inversion Problem}

Recall that the Jacobi Inversion Problem is, given a point $\bm{u} \in J$, to find the preimage of this point under the Abel map (\ref{Abel}).  

\begin{theorem} \label{JIP_Thm}
Suppose we are given $\{u_1, \dots, u_6\} = \bm{u} \in J$.   Then we could solve the Jacobi Inversion Problem explicitly using the equations derived from (\ref{pp1})-(\ref{pp5}).
\end{theorem}
\begin{pro}
Consider either equation (T1) or (T2) defined in the discussion above.  This was a polynomial constructed from $\wp$-functions and the variable $z$.  This equation had degree 6 in $z$ so denote by $(z_1,\dots,z_6)$ the 6 zeros of the polynomial.

Next, rearrange (\ref{pp1}) to give an equation for $w^2$.  Substitute this into equation (\ref{pp2}) and multiply all terms by $\wp_{66}$ to give the following equation of degree 1 with respect to $w$:
\begin{align}
0 &= w\big( z\wp_{66}\wp_{55} - 2z^2\wp_{66} - z\wp_{66}\wp_{566} + \wp_{666}z\wp_{56} + \wp_{36}\wp_{666} +
\wp_{66}\wp_{35} \nonumber \\ 
&\quad - z\wp_{56}^2 - \wp_{36}\wp_{56} - \wp_{66}\wp_{366} \big) + z^2\wp_{66}\wp_{45} - z^2\wp_{66}\wp_{466}
 + \wp_{56}z^3  \nonumber \\ 
&\quad - z\wp_{66}\wp_{266} - \wp_{56}z\wp_{26} + z\wp_{66}\wp_{25} + \wp_{15}\wp_{66} - \wp_{666}z^3 + \wp_{666}\wp_{46}z^2 \nonumber \\
&\quad + \wp_{666}z\wp_{26} - \wp_{166}\wp_{66} + \wp_{666}\wp_{16} - \wp_{56}\wp_{46}z^2 - \wp_{56}\wp_{16}. \label{ow1}
\end{align}
We could substitute each $z_i$ into equation (\ref{ow1}) in turn, and solve to find the corresponding $w_i$.  
We can therefore identify the set of points $\{(z_1, w_1), \dots, (z_6,w_6)\}$ on the curve $C$ which are the Abel preimage of $\bm{u}$. \\
\end{pro}

\section{Deriving the properties of $\bm{\sigma(u)}$}\label{SECsigexp}

In this section we derive some properties for $\sigma(\bu)$ and use them to construct the series expansion.

\begin{lemma}\label{sigzero}
The function $\sigma(\bu)$ has zeroes of order 1 when $\bu \in \Theta^{[5]}$.  Further, $\sigma(\bu) \neq 0$ for all other $\bu$.
\end{lemma}
\begin{pro}
This is a classical result, first discussed in \cite{ba97}, which always holds on $\Theta^{[g-1]}$.
The first part can also be concluded explicitly from the results of the previous section.  In Theorem \ref{JIP_Thm} we discussed how the six roots of the polynomial $(T2)$, gave us the Abel preimage of $\bu \in J$.  Now, suppose that $\bu$ is approaching $\Theta^{[5]}$, implying one of these roots is approaching infinity.  We explicitly calculate the denominator of $(T2)$ to be $\sigma(\bu)^{16}$, using Definition \ref{nordp}.  Therefore, we can conclude that when $\bu$ descends to $\Theta^{[5]}$ we must have $\sigma(\bu)=0$.  \\
\end{pro}

\noindent Consider $\bu \in \Theta^{[5]}$ which, by definition, we can express using points $P_k$ on $C$ as
\[
\bu = \int_{\infty}^{P_1} \bm{du} + \dots + \int_{\infty}^{P_5} \bm{du}.
\]
Use equations (\ref{diffs_xi}) to express $\bu$ with five local parameters.
\begin{align}
0 &= u_1 + \frac{1}{11}\xi_1^{11} + \dots + \frac{1}{11}\xi_5^{11} + O(\xi_1^{15}) + \dots + O(\xi_5^{15}), \nonumber \\
&\hspace*{0.08in} \vdots  \nonumber \\
0 &= u_6 + \xi_1 + \dots + \xi_5 + O(\xi_1^{5}) + \dots + O(\xi_5^{5}). \label{res_eq}
\end{align}
Suppose we were to take the multivariate resultant of these six equation, eliminating the parameters $\xi_1,\dots,\xi_5$.  From the theory of resultants we would be left with the unique, (up to a constant), polynomial that must be zero for $\bu \in \Theta^{[5]}$.  By Lemma \ref{sigzero} this polynomial would be equal to $\sigma(\bu)$.  Further, since it was generated using polynomials of homogeneous weight, we know that $\sigma(\bu)$ must also have definite weight.

We can perform this calculation explicitly in the case when $\bm{\lambda} = \bm{0}$.  In this case equations (\ref{diffs_xi}) simplify to 
\[
\begin{array}{llllll}
u_1 = -\frac{1}{11}\xi^{11}  & u_3= -\frac{1}{6}\xi^6  & u_5 = -\frac{1}{2}\xi^2  &
u_2 = -\frac{1}{7}\xi^{7}    & u_4= -\frac{1}{3}\xi^3  & u_6 = -\xi, 
\end{array} 
\]
and hence equations (\ref{res_eq}) become finite polynomials.  

We use a multivariate resultant calculation to eliminate $\xi_1,\dots,\xi_5$ and leave the polynomial below equal to zero.
\begin{align}
&SW_{4,5} =  \textstyle \frac{1}{8382528}u_{6}^{15} + \frac{1}{336}u_{6}^{8}u_{5}^{2}u_{4} 
- \frac{1}{12}u_{6}^{4}u_{1} - \frac{1}{126}u_{6}^{7}u_{3}u_{5} - \frac{1}{6}u_{4}u_{3}u_{5}u_{6}^{4} \nonumber \\
&\quad \textstyle - \frac{1}{72}u_{4}^{3}u_{6}^{6} - \frac{1}{33264}u_{6}^{11}u_{5}^{2} 
+ \frac{1}{27}u_{5}^{6}u_{6}^{3} + \frac{2}{3}u_{4}u_{5}^{3}u_{3} - 2u_{4}^{2}u_{6}u_{3}u_{5} - u_{2}^{2}u_{6} \nonumber \\
&\quad \textstyle - \frac{2}{9}u_{5}^{3}u_{3}u_{6}^{3} - u_{4}u_{3}^{2} + \frac{1}{12}u_{4}^{4}u_{6}^{3} 
- \frac{1}{3024}u_{6}^{9}u_{4}^{2} - \frac{1}{756}u_{6}^{7}u_{5}^{4} + \frac{1}{1008}u_{6}^{8}u_{2} \nonumber \\
&\quad \textstyle + \frac{1}{3}u_{5}^{4}u_{2} + \frac{1}{3}u_{6}^{3}u_{3}^{2} - \frac{1}{9}u_{4}u_{5}^{6} 
+ \frac{1}{399168}u_{6}^{12}u_{4} + u_{4}u_{6}u_{5}^{2}u_{2} + \frac{1}{4}u_{4}^{5} \nonumber \\
&\quad \textstyle + 2\,u_{{5}}u_{{3}}u_{{2}} + \frac{1}{6}\,{u_{{5}}}^{2}{u_{{6}}}^{4}u_{{2}} 
+ \frac{1}{12}\,{u_{{6}}}^{5}u_{{2}}u_{{4}} - \frac{1}{2}\,{u_{{4}}}^{2}{u_{{6}}}^{2}u_{{2}} 
+ \frac{1}{2}\,{u_{{4}}}^{3}{u_{{6}}}^{2}{u_{{5}}}^{2} \nonumber \\  
&\quad \textstyle - \frac{1}{3}\,{u_{{4}}}^{2}u_{{6}}{u_{{5}}}^{4}  - \frac{1}{36}\,{u_{{5}}}^{4}u_{{4}}{u_{{6}}}^{4} 
+ u_{{4}}u_{{6}}u_{{1}} - {u_{{5}}}^{2}u_{{1}}. \label{SW45}
\end{align}
In fact this is just a specific case of the following result for the $\sigma$-function.
\begin{lemma} \label{SW}
Define the canonical limit of the sigma function as the value of $\sigma(\bu)$ in the case when all the curve constants are zero.  In this case the series expansion of $\sigma(\bu)$ about $\bu = (0,0,0,0,0,0)$ is given by a constant $K$ multiplied by the Schur-Weierstrass polynomial generated by $(n,s)$.
\end{lemma}
\begin{pro}
The result was first stated in \cite{bel99}, with an alternative proof now available in \cite{N08}. \\
\end{pro}

\noindent Note that calculating $SW_{4,5}$ as the Schur-Weierstrass polynomial is, computationally, far easier than using a multivariate resultant method. 
\begin{corollary} \label{sig_odd}
The function $\sigma(\bu)$ associated with the (4,5)-curve, is odd with respect to $\bu \mapsto [-1]\bu$.
\end{corollary}
\begin{pro}
Fix the matrix $M$ that satisfies (\ref{Legendre}) and (\ref{imomeg}).  Then the solutions to (\ref{quas}) form a one dimensional space over $\C$, (see \cite{bego08} p456).  

Since we can express $\sigma([-1]\bu)$ as 
$\sigma(\bu + \ell)$ for some $\ell \in \Lambda$, it follows that both $\sigma(\bu)$ and $\sigma([-1]\bu)$ satisfy (\ref{quas}).  Therefore we have $\sigma([-1]\bu) = k\sigma(\bu)$, for some $k \in \C$.  
If we let $\bu = [-1]\bu$ then we see have $k^2 = \pm 1$, and so $\sigma(\bu)$ is either odd or even with respect to $\bu \mapsto [-1]\bu$.

We can easily check that $SW_{4,5}$ is an odd polynomial from equation (\ref{SW45}), and therefore by Lemma \ref{SW} we conclude that $\sigma(\bu)$ is an odd function. \\
\end{pro}

We now aim to derive a Taylor series expansion for $\sigma(\bu)$, similar to that of the elliptic case in equation (\ref{elliptic_sigexp}).  This will depend on the variables $\bm{u} = (u_1,\dots,u_6)$ and the curve constants $\bm{\lambda}=(\lambda_4,\dots,\lambda_0)$.  From Lemma \ref{SW} we already have, 
\begin{align*}
\sigma(\bu) = K \cdot SW_{4,5} 
  \quad + \quad \mbox{terms with degree in } \{\lambda_4, \dots \lambda_0\} \mbox{ greater than zero},
\end{align*}
for some constant $K$.  Further, we know by Corollary \ref{sig_odd} that the expansion will be odd, and also, since we know that $\sigma(\bu)$ has definite weight, we have that the expansion is homogeneous in the Sato weights.  From equation (\ref{SW}) we can see that weight is $+15$.
\begin{remark}\label{remfixc}
For simplicity, we fix the constant $c$ in Definition (\ref{sigdef}) to be the constant that makes $K=1$ in Lemma (\ref{SW}).  Note that some other authors working in this area would define $c$ to be 
\[
c = \Big( \frac{\pi^6}{\mbox{det}(w')} \Big)^{\frac{1}{2}} \cdot \frac{1}{\sqrt[4]{D}},
\]
where $D$ is the discriminant of the curve $C$.  Note that this constant cancels in the definitions of all the Abelian functions defined from $\sigma(\bu)$. Hence, any relations between such functions are independent of the choice of $c$.
\end{remark}

We now have enough information to define the following expansion for $\sigma(\bu)$.
\begin{theorem}
The function $\sigma(\bu)$ associated with (\ref{C}) has an expansion of the following form.
\[
\sigma(\textbf{u}) = \sigma(u_1, u_2, u_3, u_4, u_5, u_6) = C_{15}(\bu) +  C_{19}(\bu) +  \dots + C_{15 + 4n}(\bu) + \dots
\]
where each $C_k$ is a finite, odd polynomial composed of products of monomials in $u_i$ of total weight $+k$, multiplied by monomials in $\lambda_j$ of total weight $15-k$.
\end{theorem}
\begin{pro}
The theoretical part of the proof follows \cite{eemop07} and \cite{bego08}.   The rational is that although the expansion is homogeneous of weight $+15$, it will contain both $u_i$ (with $+$ve weight) and $\lambda_j$ (with -ve weight).  We hence split up the infinite expansion, into finite polynomials whose terms share common weight ratios.

The first polynomial will be the terms with the lowest weight in $u_i$.  These must be the terms that do not vary with $\bm{\lambda}$.  The indices then increase by four since the weights of $\bm{\lambda}$ decrease by four (see Definition \ref{Sato1}). \\
\end{pro}

By Lemma \ref{SW} we have $C_{15} = SW_{4,5}$ as given by equation (\ref{SW45}).  Using the computer algebra package Maple, we calculate the other polynomials successively as follows:
\begin{enumerate}
\item Select the terms that could appear in $C_k$.  These are a finite number of monomials formed by entries of $\bu$ and $\bm{\lambda}$ with the appropriate weight ratio.
\item Construct $\hat{\sigma}(\bu)$ as the sum of $C_k$ derived thus far.  Then add to this each of the possible terms, multiplied by an independent, unidentified constant.
\item Determine the constants by ensuring $\hat{\sigma}(\bu)$ satisfies known properties of the $\sigma$-function.
\begin{itemize}
\item For the first few $C_k$ this was mainly ensuring Lemma \ref{sigzero} is satisfied (as in the trigonal calculations).
\item For the latter $C_k$, the coefficients were instead determined by ensuring a variety of the equations from Lemma \ref{Qcor} were satisfied. 
\item In addition, those polynomials up to $C_{39}$ required we ensure $\sigma(\bu)$ satisfied some of the relations between $\wp$-functions obtained from the expansion of the Kleinian formula in Section \ref{SECkleinform}.
\end{itemize}
The second method is the most computationally efficient (due to the pole cancellations), while the third method is extremely difficult.  The equations in Lemma \ref{Qcor} are derived in tandem with the $\sigma$-function expansion, and so cannot be used for the first few $C_k$.
\end{enumerate}
The expansion has been calculated up to and including $C_{59}$.  Appendix \ref{APPsig} contains $C_{19}$ and $C_{23}$ with the rest of the expansion online at \cite{myweb}.  These latter polynomials are extremely large, and represent a significant amount of computation.  Many of the calculations were run in parallel on a cluster of machines using the Distributed Maple package (see \cite{DM}).  This expansion is sufficient for any explicit calculations.  However, it would be ideal to find a recursive construction of the expansion generalising the elliptic case, (see for example \cite{ee99}).

\section{Relations between the Abelian functions}\label{SECbasis}

In the previous section we showed that $\sigma(\bu)$ has definite Sato weight, and hence so does the Abelian functions defined from it.  We can conclude from Definition \ref{nordp} that 
\begin{equation} \label{pwt}
\mbox{wt}(\wp_{i_1,i_2,\dots,i_n}) = - \left[{\text wt}(u_{i_1})+{\text wt}(u_{i_2})+\dots{\text wt}(u_{i_n})\right].
\end{equation}
Then use equations (\ref{4iQ}) and (\ref{6iQ}) respectively to conclude 
\begin{equation}
\mbox{wt}(Q_{ijkl}) = \mbox{wt}(\wp_{ijkl}) \qquad \mbox{and} \qquad \mbox{wt}(Q_{ijklmn}) = \mbox{wt}(\wp_{ijklmn}).
\end{equation}

We now introduce the following definition to classify the Abelian functions associated with $C$ by their pole structure.
\begin{definition} Define
\[
\Gamma \big( J, \mathcal{O}(m \Theta^{[k]} ) \big)
\]
as the vector space of Abelian functions defined upon $J$, which have poles of order at most $m$, occurring only on the $k$th standard theta subset, $\Theta^{[k]}$.
\end{definition}
Recall that the Abelian functions we define all had poles occurring only when $\sigma(\bu)=0$, which by Lemma \ref{sigzero}, is when $\bu \in \Theta^{[5]}$.  Therefore, using Remark \ref{Qspec}(iii), we conclude that the $n$-index $\wp$-functions belong to $\Gamma \big( J, \mathcal{O}(n \Theta^{[5]} ) \big)$, while the $n$-index $Q$-functions all belong to $\big( J, \mathcal{O}(2 \Theta^{[5]} ) \big)$.

\begin{theorem} \label{basiseq}
A basis for $\Gamma \big( J, \mathcal{O}(2 \Theta^{[2]} ) \big)$ is given by
\begin{align*} 
\begin{array}{ccccccccccccc}
    & \C1         &\op& \C\wp_{11}  &\op& \C\wp_{12}  &\op& \C\wp_{13}  &\op& \C\wp_{14}  \\
\op & \C\wp_{15}  &\op& \C\wp_{16}  &\op& \C\wp_{22}  &\op& \C\wp_{23}  &\op& \C\wp_{24}  \\
\op & \C\wp_{25}  &\op& \C\wp_{26}  &\op& \C\wp_{33}  &\op& \C\wp_{34}  &\op& \C\wp_{35}  \\
\op & \C\wp_{36}  &\op& \C\wp_{44}  &\op& \C\wp_{45}  &\op& \C\wp_{46}  &\op& \C\wp_{55}  \\
\op & \C\wp_{56}  &\op& \C\wp_{66}  &\op& \C Q_{5566} &\op& \C Q_{4556} &\op& \C Q_{4555} \\
\op & \C Q_{4455} &\op& \C Q_{3566} &\op& \C Q_{3556} &\op& \C Q_{2566} &\op& \C Q_{2556} \\
\end{array}
\end{align*}
\begin{align*}
\begin{array}{ccccccccccccc}
\op & \C Q_{3456} &\op& \C Q_{2456} &\op& \C Q_{3366} &\op& \C Q_{3445} &\op& \C Q_{2366} \\
\op & \C Q_{2445} &\op& \C Q_{1466} &\op& \C Q_{1556} &\op& \C Q_{2266} &\op& \C Q_{2356} \\
\op & \C Q_{2256} &\op& \C Q_{2346} &\op& \C Q_{1455} &\op& \C Q_{2345} &\op& \C Q_{3344} \\
\op & \C Q_{2245} &\op& \C Q_{2344} &\op& \C Q_{1266} &\op& \C Q_{1356} &\op& \C Q_{1444} \\
\op & \C Q_{1346} &\op& \C Q_{2236} &\op& \C Q_{2335} &\op& \C Q_{1246} &\op& \C Q_{1255} \\
\op & \C Q_{1245} &\op& \C Q_{1166} &\op& \C Q_{1244} &\op& \C Q_{1156} &\op& \C Q_{1146} \\
\op & \C Q_{1155} &\op& \C Q_{1145} &\op& \C Q_{1144} &\op& \multicolumn{2}{l}{\C Q_{114466}.}
\end{array}
\end{align*}
\end{theorem}
\begin{pro}
The dimension of the space is $2^g=2^6=64$ by the Riemann-Roch theorem for Abelian varieties.  It was shown above that all the selected elements do in fact belong to the space.  All that remains is to prove their linear independence, which can be done explicitly using Maple. \\
\end{pro}

The actual construction of the basis was as follows.  We started by including all 21 of the $\wp_{ij}$ in the basis, since they were all linearly independent.  Then, to decide which $Q_{ijkl}$ to include, we systematically considered decreasing weights in turn, starting at $-4$ since this is the highest weight of any $Q$-function.  At each stage we derived equations to express the $Q$-functions at that weight using the following method (implemented with Maple):
\begin{enumerate}
\item Choose the possible terms at this weight.  These are the elements currently in the basis with this weight, along with elements in the basis of a higher weight (already determined) combined with appropriate $\lambda$-monomials, that balance the weight.  
\item We form a sum of these terms, each multiplied by an undetermined coefficient.  We also include in this sum, the $Q_{ijkl}$ which are at this weight.
\item Substitute the Abelian functions for their definitions as $\s$-derivatives.
\item Substitute $\sigma(\bu)$ for the expansion, truncated at the appropriate point.
\item Take the numerator of the resulting expression and separate into monomials in $\bm{u}$ and $\bm{\lambda}$, with coefficients in the unidentified coefficients.
\item Set all the coefficients to zero, and solve the resulting system of equations.
\end{enumerate}
At weights which have more than one $Q$-function, we often find that one or more must be added to the basis, so that the others can be expressed.

We form these equations at successively lower weights, constructing the basis as we proceed.  As the weight decreases, we require more of the expansion, which is why these were calculated in tandem.  Also, as the weight decreases the possible number of terms increase, and the computations take more time and memory.  Upon completing this process we have 63 basis elements.  

We find the final element by considering the 6-index Q-functions.  Repeating the process, we found that one of the functions at weight $-30$ is required to express the others.

\subsubsection*{Sets of differential equations satisfied by the Abelian functions}

We now present a number of differential equations between the Abelian functions.  The number in brackets on the left indicates the weight of the equation.
\begin{lemma}\label{Qcor}
Those 4-index $Q$-functions not in the basis, can be expressed as a linear combination of the basis elements.
\begin{align*}
\begin{array}{lcl}
\textbf{(-4)} \quad Q_{6666} = - 3\wp_{55} + 4\wp_{46},                                           &\hspace*{0.1in}& 
\textbf{(-7)} \quad Q_{4566} = 2\lambda_{4}\wp_{56} + 2\wp_{36},     \\
\textbf{(-5)} \quad Q_{5666} = - 2\wp_{45},                        \textcolor{white}{\vdots}      & &
\textbf{(-7)} \quad Q_{5556} = 4\lambda_{4}\wp_{56} + 4\wp_{36},     \\
\textbf{(-6)} \quad Q_{4666} = \textstyle 6\lambda_{4}\wp_{66} - 2\wp_{44} - \frac{3}{2}Q_{5566}, & & 
\hspace*{0.82in} \vdots
\end{array}
\end{align*}
A longer list is given in Appendix \ref{APPQ}, while the full set is available online at \cite{myweb}.

The same statement is also true for all the 6-index $Q$-functions, except $Q_{114466}$.  Explicit relations have been calculated down to weight $-30$.  The first few are given below, with all available relations online at \cite{myweb}.
\begin{align*}
\textbf{(-6)} \quad Q_{666666} &= 40\wp_{44} + 15Q_{5566} - 24\wp_{66}\lambda_4, \nonumber \\
\textbf{(-7)} \quad Q_{566666} &= 20\wp_{36} - 4\wp_{56}\lambda_4, \nonumber \\
\textbf{(-8)} \quad Q_{556666} &= 24\wp_{26} - 12\wp_{35} - 2Q_{4556}, \\
\textbf{(-8)} \quad Q_{466666} &\underset{\vdots}{=} - 20\wp_{35} + 5Q_{4556} + 16\wp_{46}\lambda_4 - 20\wp_{55}\lambda_4 
- 8\lambda_3, \nonumber
\end{align*}
\end{lemma}
\begin{pro}
Clearly such relations must exist.  The explicit PDEs were calculated in the construction of the basis, as discussed at the start of this section. \\
\end{pro}

\begin{corollary}\label{Pcor}
There are a set of PDEs that express 4-index $\wp$-functions, using Abelian functions of order at most 2.  The full set can be found online at \cite{myweb}.
\begin{align}
\textbf{(-4)} \quad \wp_{6666} &= 6\wp_{66}^2 - 3\wp_{55} + 4\wp_{46} \label{p6666}, \\
\textbf{(-5)} \quad \wp_{5666} &= 6\wp_{56}\wp_{66} - 2\wp_{45} \label{p5666}, \\
\textbf{(-6)} \quad \wp_{4666} &= 6\wp_{46}\wp_{66} + 6\lambda_{4}\wp_{66} - 2\wp_{44} - \textstyle \frac{3}{2}\wp_{5566} + 3\wp_{66}\wp_{55} + 6\wp_{56}^2, \nonumber \\
\textbf{(-7)} \quad \wp_{4566} &= 2\wp_{45}\wp_{66} + 4\wp_{46}\wp_{56} + 2\lambda_{4}\wp_{56} + 2\wp_{36}, \nonumber \\
\textbf{(-7)} \quad \wp_{5556} &\underset{\vdots}{=} 6\wp_{55}\wp_{56} + 4\lambda_{4}\wp_{56} + 4\wp_{36},  \nonumber
\end{align}
\end{corollary}
\begin{pro}
Apply equation (\ref{4iQ}) to the first set of relations in Lemma \ref{Qcor}. \\
\end{pro}

The set of equations in Corollary \ref{Pcor} is of particular interest because it gives a generalisation of equation (\ref{elliptic_diff1}), from the elliptic case. A similar generalisation for equation (\ref{elliptic_diff2}) would be a set of equations that express the 3-index $\wp$-functions, using Abelian function of order at most 3.  
So far the following relations have been derived (see \cite{myweb} for the latest list).
\begin{align*}
&\textbf{(-6)} \hspace*{0.4in} \wp_{666}^{2} = \textstyle 4\wp_{66}^{3} - 7\wp_{56}^{2} + 4\wp_{46}\wp_{66} 
- 8\wp_{55}\wp_{66} - 4\wp_{66}\lambda_{4} + 4\wp_{44} + 2\wp_{5566}, \\
&\textbf{(-7)} \quad \wp_{566}\wp_{666} = \textstyle 4\wp_{66}^{2}\wp_{56} + 2\wp_{46}\wp_{56} - \wp_{55}\wp_{56} 
- 2\wp_{45}\wp_{66} + 2\wp_{36}, \\
&\textbf{(-8)} \quad \wp_{556}\wp_{666} = \textstyle - 4\wp_{26} - 2\wp_{35} - 4\wp_{55}\lambda_{4} - 4\lambda_{3} 
+ 2\wp_{4556} - 6\wp_{45}\wp_{56} - 2\wp_{46}\wp_{55} \\
&\hspace*{0.7in} \textstyle + \wp_{5566}\wp_{66} - 2\wp_{56}^{2}\wp_{66} - 2\wp_{55}^{2}, \\
&\textbf{(-8)} \hspace*{0.4in} \wp_{566}^2 = \textstyle 4\wp_{56}^{2}\wp_{66} + 4\wp_{46}\wp_{55} + \wp_{55}^{2} 
+ 4\wp_{55}\lambda_{4} + 4\wp_{45}\wp_{56} +8\wp_{26} \\ 
&\hspace*{0.7in} \textstyle +4\lambda_{{3}} - 2\wp_{4556}, \\
&\textbf{(-8)} \quad \wp_{466}\wp_{666} = \textstyle 4\wp_{56}^{2}\wp_{66} + 4\wp_{46}\wp_{66}^{2} 
+ 2\wp_{55}\wp_{66}^{2} + 4\wp_{66}^{2}\lambda_{4} + 2\wp_{46}^{2} - 4\wp_{46}\wp_{55} \\
&\hspace*{0.7in} \textstyle - 2\wp_{55}\lambda_{{4}} - 3\wp_{45}\wp_{56} - 2\wp_{44}\wp_{66} - \wp_{5566}\wp_{66} 
- 2\wp_{26} - 2\wp_{35} \\ 
&\hspace*{0.7in} \textstyle - 2\lambda_{{3}}+\wp_{4556}, \\
&\textbf{(-9)} \quad \wp_{556}\wp_{566} = \textstyle -2\wp_{56}^{3} - 2\wp_{45}\wp_{55} 
- \frac{4}{3}\wp_{45}\lambda_{4} + \wp_{5566}\wp_{56} - \frac{4}{3}\wp_{34} + \frac{1}{3}\wp_{4555}, \\
&\textbf{(-9)} \quad \wp_{555}\wp_{666} = \textstyle -3\wp_{5566}\wp_{56} - 4\wp_{44}\wp_{56} 
+ 8\wp_{56}\wp_{66}\lambda_{{4}} + 12\wp_{56}\wp_{55}\wp_{66} + 10\wp_{56}^{3} \\
&\hspace*{0.7in} \textstyle + \frac{8}{3}\wp_{45}\lambda_{4} 
- \frac{4}{3}\wp_{34} - \frac{2}{3}\wp_{4555} + 4\wp_{45}\wp_{46} + 8\wp_{36} \wp_{66}, \\
&\textbf{(-9)} \quad \wp_{466}\wp_{566} = \textstyle -2\wp_{25} + \frac{2}{3}\wp_{34} - \frac{4}{3}\wp_{45}\lambda_{4} 
- 2\wp_{45}\wp_{46} + \frac{1}{3}\wp_{4555} + 2\wp_{56}\wp_{55} \wp_{66} \\
&\hspace*{0.7in} \textstyle + 4\wp_{56}\wp_{46}\wp_{66} + 4\wp_{56}^{3} 
- \wp_{5566}\wp_{56} - 2\wp_{44}\wp_{56} + 4\wp_{56}\wp_{66}\lambda_{{4}} - \wp_{45} \wp_{55}, \\
&\textbf{(-9)} \quad \wp_{456}\wp_{666} = \textstyle -\wp_{56}\wp_{55}\wp_{66} + 2\wp_{45}\wp_{66}^{2} 
+ 2\wp_{56}\wp_{46}\wp_{66} + 2\wp_{36}\wp_{66} - 2\wp_{56}^{3} \\
&\hspace*{0.7in} \textstyle + 2\wp_{44}\wp_{56} + \frac{4}{3}\wp_{45}\lambda_{{4}} 
+ \frac{1}{2}\wp_{5566}\wp_{56} - \frac{2}{3}\wp_{34} - \frac{1}{3}\wp_{4555} + 2\wp_{45}\wp_{46}. \\
\end{align*}

\begin{proposition}  There are a set of relations that are bi-linear in the 2-index and 3-index $\wp$-functions.  (See \cite{myweb} for full list.) There is no analog in the elliptic case, although similar relations have been derived in the hyperelliptic and trigonal cases. 
\begin{align}
\textbf{(-6)} \quad 0 &= -\wp_{555} + 2\wp_{456}  +2\wp_{566}\wp_{66} - 2\wp_{56}\wp_{666}, 
\label{lin06}  \\ 
\textbf{(-7)} \quad 0 &= -2\wp_{446} + 2\wp_{455} -2\wp_{466} \wp_{66} +2\wp_{666} \lambda_{4} + 2\wp_{46}\wp_{666} 
\nonumber \\
&\qquad - 2\wp_{556}\wp_{66} + \wp_{55} \wp_{666} +\wp_{566} \wp_{56}, 
 \nonumber  \\ 
\textbf{(-8)} \quad 0 &= - 2\wp_{56}\wp_{466} + 2\wp_{46}\wp_{566} + \wp_{555}\wp_{66} - 2\wp_{55}\wp_{566} 
+ \wp_{556}\wp_{56} - 2\wp_{366}, 
\nonumber  \\ 
\textbf{(-8)} \quad 0 &= - \wp_{445} +2\wp_{456} \wp_{66} +\wp_{56} \wp_{466} -\wp_{366} - \wp_{566}\lambda_{4} 
- \wp_{45}\wp_{666} - 2\wp_{46} \wp_{566}, \nonumber \\
\textbf{(-9)} \quad 0 &= - 2\wp_{455} \wp_{66} + 4\wp_{266} + 2\wp_{45} \wp_{566} + 2\wp_{466}\wp_{55} 
- 2\wp_{46}\wp_{556} + \wp_{556}\wp_{55} \nonumber \\ 
&\qquad - \wp_{555}\wp_{56} - 2\wp_{356}. \nonumber
\end{align}
\end{proposition}
\begin{pro}
These can be calculated by cross differentiating suitable pairs of equations from Corollary \ref{Pcor}.  For example,  equation (\ref{p6666}) expresses $\wp_{6666}$ while equation (\ref{p5666}) expresses $\wp_{5666}$.  If we substitute for these equations into 
\[
\frac{\partial}{\partial u_5} \wp_{6666} - \frac{\partial}{\partial u_6} \wp_{5666} = 0,
\]
then we find equation (\ref{lin06}). \\
\end{pro}

\noindent A topic of future work in this area would be the construction of relations between the $\wp$-functions in covariant form, as was recently acheived in \cite{CA2008} for the hyperelliptic case.

\section{Solution to the KP-equation} \label{SECnonlin}

We now demonstrate how such Abelian functions can give a solution to the KP-equation.  Differentiate equation (\ref{p6666}) twice with respect to $u_6$ to obtain
\begin{align*}
\wp_{666666} &= 12\tfrac{\partial}{\partial u_6}\big(\wp_{66}\wp_{666}\big) - 3\wp_{5566} + 4\wp_{4666}.
\end{align*} 
Let $u_6=x, u_5=y, u_4=t$ and $W(x,y,t) = \wp_{66}(\bm{u})$.  We then rearrange to give
\[
\big[ W_{xxx} - 12WW_x - 4W_t \big]_x + 3W_{yy} = 0,
\]
a parametrised form of the KP-equation.  In fact, this is just a special case of the following general result for Abelian functions associated with algebraic curves.

\begin{theorem}
Let $E$ be an $(n,s)$-curve with genus $g$ as given by equation (\ref{ns_curve}).   Define the multivariate $\sigma$-function associated with $E$ as normal.  Define the Abelian functions from $\sigma(\bu)$ as in Section \ref{SECabeliandef} (with the indicies now running to $g$ instead of 6). Finally, define the function $W(\bm{u}) = \wp_{gg}(\bm{u})$, which we denote $W(x,y,t)$ after applying the substitutions 
$u_g=x, u_{g-1}=y, u_{g-2}=t$.  

Then, if $n\geq4$ the function $W(x,y,t)$ will satisfy the following parametrised version of the KP-equation.
\begin{equation} \label{para_KP}
\big( W_{xxx} - 12WW_x - bW_t \big)_x - aW_{yy} = 0,
\end{equation}
for some constants $a,b$.
\end{theorem}
\begin{pro}
Recall Definitions \ref{Sato1} and \ref{Sato2} which gave the Sato weights for $C$.  These can be calculated for the general curve $E$ similarly as
\begin{align*}
\begin{array}{lllll}
\mbox{wt}(x) = -n,     & \quad & \mbox{wt}(y) = -s,        \\
\mbox{wt}(u_g) = \w_1, & \quad & \mbox{wt}(u_{g-1}) = \w_2, & \dots & \mbox{wt}(u_1) = \w_g, \\
\mbox{wt}(\Y_0) = -ns, & \quad & \mbox{wt}(\Y_1) = -n(s-1), & \dots & \mbox{wt}(\Y_{s-1}) = -n.
\end{array}
\end{align*}
Here $n,s$ are the integers generating the curve, and $\{\w_1,\dots,\w_g\}$ is the Weierstrass gap sequence for $n,s$.  These are the natural numbers not representable in the form $an + bs$ where $a,b \in \N$. (See \cite{bel99} Section 1 for more details.)

Since $s>n>4$ we know that $\{1,2,3\}$ cannot be represented in this form.  Therefore, we have
\[
\mbox{wt}(u_g) = +1, \qquad  \mbox{wt}(u_{g-1}) = +2, \qquad \mbox{wt}(u_{g-2}) = +3.
\]
By equation (\ref{pwt}), this implies the $\wp$-functions will have weights
\[
\mbox{wt}(\wp_{g,g}) = -2, \quad \mbox{wt}(\wp_{g,g-1}) = -3, \quad \mbox{wt}(\wp_{g-1,g-1}) = -4, \quad \mbox{wt}(\wp_{g-2,g}) = -4,
\]
with all the other 2-index $\wp$-functions having a lower weight.  
Next consider $Q_{gggg}$, which will have weight $-4$.  
This will belong to, $\Gamma \big( J, \mathcal{O}(2 \Theta^{[g-1]} ) \big)$, the space of Abelian functions defined upon the Jacobian of $E$, which have poles of at most order $2$ on $\Theta^{[g-1]}$.  We can therefore express $Q_{gggg}$ as
\begin{align*}
Q_{gggg} = a\wp_{g-1,g-1} + b\wp_{g-2,g}, \qquad (a,b \mbox{ constants}),
\end{align*}
since these are the only Abelian functions of weight $-4$.  We use Remark \ref{Qspec} to substitute for $Q$, and then differentiate twice with respect to $u_g$ to give
\begin{align*}
\wp_{gggggg} &= 12\frac{\partial}{\partial u_g}\big(\wp_{gg}\wp_{ggg}\big) + a\wp_{g-1,g-1,g,g} + b\wp_{g-2,g,g,g}.
\end{align*}
Then make the substitutions suggested in the theorem to obtain equation (\ref{para_KP}). 
\hspace*{0.1in}

\end{pro} 

Further research into the possible applications of these results is currently being conducted.

\section{Two-term addition formula}\label{SECaddform}

\begin{theorem}
The functions associated with (\ref{C}) satisfy the following two-term addition formula:
\[
- \frac{\s(\bm{u}+\bm{v})\s(\bm{u}-\bm{v})}{\s(\bm{u})^2\s(\bm{v})^2} = f(\bm{u},\bm{v}) - f(\bm{v},\bm{u}),
\]
where $f(\bu,\bv)$ is a polynomial of Abelian functions, given in Appendix \ref{APPadd}
\end{theorem}
\begin{pro}
We seek to express the following ratio of sigma functions (labelled LHS), using a sum of Abelian functions.
\begin{equation}
\mbox{LHS}(\bm{u},\bm{v}) = - \frac{\s(\bm{u}+\bm{v})\s(\bm{u}-\bm{v})}{\s(\bm{u})^2\s(\bm{v})^2}.
\end{equation}
First recall that $\sigma(\bu)$ is an odd function, with respect to the change of variables $\bm{u} \mapsto [-1]\bm{u}$.  We use this to consider the effect of $(\bm{u},\bm{v}) \mapsto (\bm{v},\bm{u})$ on LHS.
\begin{align*}
\mbox{LHS}(\bm{v},\bm{u}) &= - \frac{\s(\bm{v}+\bm{u})\s(\bm{v}-\bm{u})}{\s(\bm{v})^2\s(\bm{u})^2}
                          = - \frac{\s(\bm{u}+\bm{v})\s\big([-1](\bm{u}-\bm{v})\big)}{\s(\bm{u})^2\s(\bm{v})^2}
                          = - \mbox{LHS}(\bm{u},\bm{v}).
\end{align*}
So LHS is antisymmetric, or odd with respect to $(\bm{u},\bm{v}) \mapsto (\bm{v},\bm{u})$. 

Next, recall that $\sigma(\bm{u})$ has zeros of order 1 along $\Theta^{[5]}$ and no zeros anywhere else.  This implies that LHS has poles of order 2 along 
\[
(\Theta^{[5]} \times J) \cup (\Theta^{[5]} \times J)
\]
but nowhere else.  Together, this implies that we can express $LHS$ as
\begin{equation} \label{LHSeq}
LHS = \sum_j A_j \Big( X_j(\bu)Y_j(\bv) - X_j(\bv)Y_j(\bu) \Big),
\end{equation}
where the $A_j$ are constant coefficients (which may be functions of $\bm{\lambda}$), and the $X_j$ and $Y_j$ are functions chosen from the basis in Theorem \ref{basiseq}.

Finally we use the fact that sigma has weight $+15$ to determine that the weight of LHS is $-30$.  Hence we need only consider those terms in equation (\ref{LHSeq}) that give the correct overall weight.

We use Maple to construct equation (\ref{LHSeq}) with the $A_j$ undetermined.  This contained 1348 terms (647 undetermined coefficients since it is antisymmetric).  The coefficients were determined using the $\sigma$-function expansion.\\
\end{pro}

We believe this to be the first of a family of similar addition formula, related to the invariance expressed in equation (\ref{invar}).  There has been much work conducted into these addition formula for the trigonal cases (see \cite {eemop07} for example).  In \cite{emo08} we see that this has inspired new results in the lower genus cases.

\section*{Acknowledgements}

We would like to thank Dr.~Y.~Onishi and Dr.~A.~Nakayashiki for their suggestion that assisted in the construction of the basis in Theorem \ref{basiseq}.  We would also like to that Dr.~J.~Gibbons and Dr.~V.~Enolski for useful conversations. 

\newpage

\begin{footnotesize}
\bibliography{45-Paper}{}
\bibliographystyle{plain}
\end{footnotesize}

\vspace*{0.1in}

\begin{footnotesize}
\noindent \textbf{Department of Mathematics and the Maxwell Institute for Mathematical Sciences,
Heriot-Watt University, Edinburgh, UK EH14 4AS \\
E-mail address: \textit{M.England@ma.hw.ac.uk}, \textit{J.C.Eilbeck@hw.ac.uk}}
\end{footnotesize}

\appendix

\section{Expansion of the Kleinian formula}\label{APP_KE}

We consider equation (\ref{Kleineq}) from Theorem \ref{Klein_Thm}.  We expand this as one of the $P_k$ tends to infinity, to obtain a series expansion in terms of the local parameter $\xi$.  It follows that each coefficient with respect to $\xi$ must be zero, giving us an infinite sequence of equations, starting with those below. 
\begin{align}
0 &= \textstyle - {z}^{3} + \wp_{{46}}{z}^{2} +  \left( \wp_{{56}}w + \wp_{{26}} \right) z + \wp_{{66}}{w}^{2} 
+ \wp_{{36}}w + \wp_{{16}}  \label{pp1} \\
0 &= \textstyle  \left( \wp_{{45}} - \wp_{{466}} - 2w \right) {z}^{2} 
+  \left(  \left( \wp_{{55}} - \wp_{{566}} \right) w + \wp_{{25}} - \wp_{{266}} \right) z 
+  \left( \wp_{{56}} - \wp_{{666}} \right) {w}^{2} \nonumber \\
&\quad \textstyle +  \left( \wp_{35} - \wp_{{366}} \right) w + \wp_{{15}} - \wp_{{166}} \label{pp2}  \\
0 &= \textstyle \left( \wp_{{44}} - \frac{3}{2}\wp_{{456}} + \frac{1}{2}\wp_{{4666}}  \right) {z}^{2}
+  \big( - 3{w}^{2} +  \left( \wp_{{45}} - \frac{3}{2}\wp_{{556}} + \frac{1}{2}\wp_{{5666}} \right) w 
- \frac{3}{2}\wp_{{256}} \nonumber \\
&\quad \textstyle  + \wp_{{24}} + \frac{1}{2}\wp_{{2666}} \big)z 
+  \left( \frac{1}{2}\wp_{{6666}} - \frac{3}{2}\wp_{{566}} + \wp_{{46}} \right) {w}^{2} 
+  \big( \frac{1}{2}\wp_{{3666}} + \wp_{{34}}  \nonumber \\
&\quad \textstyle - \frac{3}{2}\wp_{{356}} \big) w + \wp_{{14}} + \frac{1}{2}\wp_{{1666}} - \frac{3}{2}\wp_{{156}} 
\label{pp3} \\
0 &= \textstyle \left( \wp_{{4566}} - \frac{4}{3}\wp_{{446}} - \frac{1}{6}\wp_{{46666}} 
- \frac{1}{2}\wp_{{455}} \right) {z}^{2} +  \big( \big( \wp_{{5566}} - \frac{4}{3}\wp_{{456}}  
- \frac{1}{6}\wp_{{56666}} \nonumber \\
&\quad \textstyle - \frac{1}{2}\wp_{{555}} \big) w - \frac{4}{3}\wp_{{246}} - \frac{1}{2}\wp_{{255}} 
+ \wp_{{2566}} - \frac{1}{6}\wp_{{26666}} \big) z - 4{w}^{3} 
+ \big( \wp_{{5666}} - \frac{1}{2}\wp_{{556}} \nonumber \\
&\quad \textstyle - \frac{1}{6}\wp_{{66666}} - \frac{4}{3}\wp_{{466}} \big) {w}^{2} 
+  \left(  - \frac{1}{6}\wp_{{36666}} + \wp_{{3566}} - \frac{4}{3}\wp_{{346}} - \frac{1}{2}\wp_{{355}} \right) w 
\nonumber \\ 
&\quad \textstyle - \frac{1}{6}\wp_{{16666}} + \wp_{{1566}} - \frac{4}{3}\wp_{{146}} - \frac{1}{2}\wp_{{155}} 
\label{pp4} \\
0 &= \textstyle - 3{z}^{4} + \left( - 2\wp_{{46}} - \frac{9}{2}\lambda_{{4}} \right) {z}^{3} 
+ \big( \frac{5}{8}\wp_{{4556}} - \frac{5}{6}\wp_{{445}} - {\frac {5}{12}}\wp_{{45666}} 
- 2\wp_{{56}}w  - 2\wp_{{26}} \nonumber \\
&\quad \textstyle - 3\lambda_{{3}} + \frac{5}{6}\wp_{{4466}} 
 + \frac{1}{2}\wp_{{46}}\lambda_{{4}} + \frac{1}{24}\wp_{{466666}} \big) {z}^{2} 
+ \big( - 2\wp_{{66}}{w}^{2} +  \big( \frac{1}{24}\wp_{{566666}} \nonumber \\
&\quad \textstyle - \frac{5}{6}\wp_{{455}}  + \frac{5}{6}\wp_{{4566}} - {\frac{5}{12}}\wp_{{55666}} + \frac{1}{2}\wp_{{56}}\lambda_{{4}} + \frac{5}{8}\wp_{{5556}} - 2\wp_{{36}} \big) w
- {\frac {5}{12}}\wp_{{25666}} \nonumber \\
&\quad \textstyle - \frac{5}{6}\wp_{{245}} - 2\lambda_{{2}} + \frac{5}{6}\wp_{{2466}} + \frac{5}{8}\wp_{{2556}}  + \frac{1}{24}\wp_{{266666}} - 2\wp_{{16}} + \frac{1}{2}\wp_{{26}}\lambda_{{4}} \big) z \nonumber \\
&\quad \textstyle +  \left( \frac{1}{2}\wp_{{66}}\lambda_{{4}} - {\frac {5}{12}}\wp_{{56666}}  
+ \frac{1}{24}\wp_{{666666}} - \frac{5}{6}\wp_{{456}} + \frac{5}{6}\wp_{{4666}} 
+ \frac{5}{8}\wp_{{5566}} \right) {w}^{2} \nonumber \\
&\quad \textstyle +  \left( \frac{1}{24}\wp_{{366666}} + \frac{5}{6}\wp_{{3466}} + \frac{5}{8}\wp_{{3556}} 
+ \frac{1}{2}\wp_{{36}}\lambda_{{4}} - {\frac {5}{12}}\wp_{{35666}} - \frac{5}{6}\wp_{{345}} \right) w 
 \nonumber \\
&\quad \textstyle + \frac{5}{8}\wp_{{1556}} - \lambda_{{1}} - {\frac {5}{12}}\wp_{{15666}} + \frac{1}{24}\wp_{{166666}} + \frac{5}{6}\wp_{{1466}} - \frac{5}{6}\wp_{{145}} + \frac{1}{2}
\wp_{{16}}\lambda_{{4}} \label{pp5} 
\end{align}

\newpage

\section{The $\bm{\sigma}$-function expansion}\label{APPsig}

We defined the $\sigma$-function expansion as the infinite sum of finite polynomials
\[
\sigma(\textbf{u}) = C_{15} + C_{19} + C_{23} + C_{27} + C_{31} + C_{35} + \dots
\]
We know that $C_{15}$ was equal to the Schur-Weierstrass polynomial as given in equation (\ref{SW45}).  The other polynomials were calculated in turn using the method described in Section \ref{SECsigexp}.  The next polynomials, $C_{19}$ and $C_{23}$, are given below, while the rest of the expansion can be found at \cite{myweb}.

\begin{align*}
&C_{19} = \textstyle \lambda_4 \cdot \big[ {\frac {1}{13970880}}\,{u_{{6}}}^{16}u_{{4}} + {\frac {2}{135}}\,{u_{{6}}}^{3}{u_{{5}}}^{8}-{u_{{6}}}^{2}u_{{5}}{u_{{4}}}^{3}u_{{3}}-
{u_{{6}}}^{2}u_{{4}}{u_{{2}}}^{2} +{\frac {4}{
45}}\,{u_{{5}}}^{6}u_{{2}} \\
&\quad \textstyle - {\frac{2}{45}}\,{u_{{6}}}^{3}{u_{{5}}}^{5}u_{{
3}} - \frac{1}{5}\,u_{{6}}{u_{{5}}}^{6}{u_{{4}}}^{2}+{\frac {1}{90}}\,{u_{{6}}}^
{4}{u_{{5}}}^{6}u_{{4}}-\frac{1}{30}\,{u_{{6}}}^{5}u_{{4}}u_{{1}} +{\frac {5}{3024}}\,{u_{{6}}}^{9}{u_{{5}}}^{2}{u_{{4}}}^{2} \\
&\quad \textstyle + \frac{1}{10}\,{u_{{5}}}^{2}{u_{{4}}}^{5}+\frac{1}{30}\,{u_{
{6}}}^{5}{u_{{2}}}^{2}+{\frac {1}{20956320}}\,{u_{{6}}}^{15}{u_{{5}}}^
{2}-\frac{1}{3}\,{u_{{5}}}^{4}u_{{1}}-{\frac {2}{45}}\,{u_{{5}}}^{8}u_{{4}}+\frac{1}{2}\,{u_{{4}}}^{4}u_{{2}} \\
&\quad \textstyle -{\frac {1}{1890}}\,{u_{{6}}}^{7}{u_{{5}}}^{6}-
{\frac {1}{120}}\,{u_{{6}}}^{7}{u_{{4}}}^{4}+{\frac {1}{665280}}\,{u_{
{6}}}^{13}{u_{{4}}}^{2}-{\frac {1}{83160}}\,{u_{{6}}}^{11}{u_{{5}}}^{4
}+{\frac {1}{997920}}\,{u_{{6}}}^{12}u_{{2}} \\
&\quad \textstyle -{\frac {1}{630}}\,{u_{{6}
}}^{7}{u_{{3}}}^{2}-{\frac {1}{5040}}\,{u_{{6}}}^{10}{u_{{4}}}^{3}+{
\frac {1}{3024}}\,{u_{{6}}}^{9}u_{{4}}u_{{2}} + \frac{1}{2}\,{u_{{6}}}^{2}{u_{{4
}}}^{2}u_{{1}}+{u_{{6}}}^{2}{u_{{5}}}^{2}{u_{{4}}}^{2}u_{{2}} \\
&\quad \textstyle -{\frac{1}{280}}\,{u_{{6}}}^{8}u_{{5}}u_{{4}}u_{{3}}+\frac{2}{15}\,{u_{{6}}}^{5}{u_{{5
}}}^{2}u_{{4}}u_{{2}}-\frac{2}{9}\,{u_{{6}}}^{4}{u_{{5}}}^{3}u_{{4}}u_{{3}}-\frac{2}{3}\,u_{{6}}{u_{{5}}}^{3}{u_{{4}}}^{2}u_{{3}} \\
&\quad \textstyle -\frac{1}{10}\,{u_{{6}}}^{5}u_{{5}}
{u_{{4}}}^{2}u_{{3}}+\frac{1}{3}\,u_{{6}}{u_{{5}}}^{4}u_{{4}}u_{{2}}-u_{{6}}{u
_{{4}}}^{2}{u_{{3}}}^{2}+\frac{1}{20}\,{u_{{6}}}^{4}{u_{{4}}}^{5}+{\frac {1}{
2520}}\,{u_{{6}}}^{8}{u_{{5}}}^{4}u_{{4}} \\
&\quad \textstyle + \frac{1}{30}\,{u_{{6}}}^{6}{u_{{4}}}
^{2}u_{{2}}+\frac{1}{6}\,{u_{{6}}}^{4}u_{{4}}{u_{{3}}}^{2}+\frac{4}{3}\,{u_{{5}}}^{3}u
_{{3}}u_{{2}}-{\frac {4}{945}}\,{u_{{6}}}^{7}{u_{{5}}}^{3}u_{{3}}+\frac{1}{18}\,{u_{{6}}}^{4}{u_{{5}}}^{4}u_{{2}} \\ &\quad \textstyle -{\frac {1}{180}}\,{u_{{6}}}^{6}{u_
{{5}}}^{2}{u_{{4}}}^{3}-{\frac {17}{997920}}\,{u_{{6}}}^{12}{u_{{5}}}^
{2}u_{{4}}+{\frac {1}{630}}\,{u_{{6}}}^{8}{u_{{5}}}^{2}u_{{2}}+{\frac 
{3}{20}}\,u_{{6}}{u_{{4}}}^{6}+{\frac {1}{5040}}\,{u_{{6}}}^{8}u_{{1}} \\
&\quad \textstyle - \frac{1}{6}\,{u_{{6}}}^{3}{u_{{4}}}^{3}u_{{2}}-{\frac {1}{83160}}\,{u_{{6}}}^
{11}u_{{5}}u_{{3}}-{\frac {1}{60}}\,{u_{{6}}}^{5}{u_{{5}}}^{4}{u_{{4}}
}^{2}+ \frac{1}{3}\,{u_{{6}}}^{3}{u_{{5}}}^{2}{u_{{4}}}^{4}+\frac{2}{15}\,{u_{{5}}}^{5}
u_{{4}}u_{{3}} \big].
\end{align*}

\begin{align*}
&C_{23}(\bm{u}) = \textstyle \lambda_3 \cdot \big[  {u_{{5}}}^{2}{u_{{3}}}^{2}u_{{2}}-{\frac{1}{360}}\,{u_{{6}}}^{6}{u_{{4}}}^{2}u_{{1}} + \frac{1}{6}\,{u_{{5}}}^{4}u_{{4}}{u_{{3}}}^{2} - \frac{1}{18}\,{u_{{6}}}^{3}{u_{{5}}}^{4}{u_{{3}}}^{2} \\ 
&\quad \textstyle - {\frac{1}{252}}\,{u_{{6}}}^{7}{u_{{5}}}^{2}{u_{{3}}}^{2} 
 +{\frac{1}{72}}\,{u_{{6}}}^{4}{u_{{5}}}^{4}u_{{1}} + \frac{1}{6}\,u_{{6}}{u_{{5}}}^{4}{u_{{2}}}^{2} + \frac{1}{12}\,{u_{{6}}}^{5}{u_{{5}}}^{2}{u_{{2}}}^{2} \\
&\quad \textstyle + \frac{1}{12}\,{u_{{6}}}^{4}{u_{{3}}}^{2}u_{{2}} + {\frac{1}{2016}}\,{u_{{6}}}^{8}{u_{{5}}}^{2}u_{{1}} 
{\frac {5}{144}}\,{u_{{6}}}^{4}{u_{{4}}}^{4}u_{{2}}  -{\frac {29}{17962560}}\,{u_{{6}}}^{12}{u_{{5}}}^{4}u_{{4}} \\
&\quad \textstyle + {\frac {1}{111767040}}\,{u_{{6}}}^{16}{u_{{5}}}^{2}u_{{4}} 
+ {\frac{17}{11975040}}\,{u_{{6}}}^{13}u_{{4}}u_{{2}} 
-  {\frac {17}{315}}\,{u_{{5}}}^{7}u_{{4}}u_{{3}} 
-  \frac{1}{30}\,{u_{{6}}}^{5}u_{{5}}u_{{4}}u_{{3}}u_{{2}} \\
&\quad \textstyle + \frac{1}{18}\,{u_{{6}}}^{3}{u_{{4}}}^{3}u_{{1}} +  \frac{1}{5}\,{u_{{5}}}^{5}u_{{3}}u_{{2}} -  {\frac{1}{74844}}\,{u_{{6}}}^{12}{u_{{5}}}^{2}u_{{2}} -  \frac{1}{6}\,{u_{{6}}}^{3}{u_{{4}}}^{2}{u_{{2}}}^{2} \\
&\quad \textstyle + {\frac {1}{630}}\,u_{{6}}{u_{{5}}}^{8}{u_{{4}}}^{2} + {\frac{1}{7560}}\,{u_{{6}}}^{4}{u_{{5}}}^{8}u_{{4}} -  {\frac {1}{45360}}\,{u_{{6}}}^{8}{u_{{5}}}^{6}u_{{4}} -  {\frac{1}{60}}\,u_{{6}}{u_{{5}}}^{2}{u_{{4}}}^{6} \\ 
&\quad \textstyle + {\frac {1}{240}}\,{u_{{6}}}^{4}{u_{{5}}}^{2}{u_{{4}}}^{5} + {\frac{1}{216}}\,{u_{{6}}}^{3}{u_{{5}}}^{4}{u_{{4}}}^{4} + {\frac {1}{1680}}\,{u_{{6}}}^{7}{u_{{5}}}^{2}{u_{{4}}}^{4} + \frac{2}{15}\,u_{{5}}{u_{{4}}}^{5}u_{{3}} \\
&\quad \textstyle + {\frac{1}{6480}}\,{u_{{6}}}^{6}{u_{{5}}}^{4}{u_{{4}}}^{3} + {\frac{1}{36288}}\,{u_{{6}}}^{10}{u_{{5}}}^{2}{u_{{4}}}^{3} + {\frac {1}{299376}}\,{u_{{6}}}^{11}{u_{{5}}}^{3}u_{{3}} -  \frac{1}{2}\,{u_{{6}}}^{2}{u_{{2}}}^{3} \\
&\quad \textstyle -  {\frac {1}{540}}\,{u_{{6}}}^{2}{u_{{5}}}^{6}{u_{{4}}}^{3} + {\frac{1}{25147584}}\,{u_{{6}}}^{15}u_{{5}}u_{{3}} -  {\frac {1}{1260}}\,{u_{{6}}}^{7}{u_{{5}}}^{5}u_{{3}}
+ \frac{1}{6}\,{u_{{4}}}^{3}{u_{{2}}}^{2} \\
&\quad \textstyle + {\frac{1}{335301120}}\,{u_{{6}}}^{17}{u_{{4}}}^{2} -  {\frac{1}{630}}\,{u_{{5}}}^{8}u_{{2}} -  {\frac {1}{1890}}\,{u_{{5}}}^{10}u_{{4}} + {\frac{1}{60}}\,{u_{{6}}}^{5}{u_{{4}}}^{2}{u_{{3}}}^{2} \\
&\quad \textstyle -  {\frac {11}{10080}}\,{u_{{6}}}^{8}u_{{4}}{u_{{3}}}^{2} + {\frac {7}{60}}\,u_{{6}}{u_{{4}}}^{5}u_{{2}} -  {\frac{37}{181440}}\,{u_{{6}}}^{10}{u_{{4}}}^{2}u_{{2}} -  {\frac{13}{1890}}\,{u_{{6}}}^{7}{u_{{4}}}^{3}u_{{2}} \\
&\quad \textstyle + {\frac{17}{945}}\,{u_{{6}}}^{3}{u_{{5}}}^{7}u_{{3}} + \frac{1}{12} \,{u_{{5}}}^{2}{u_{{4}}}^{4}u_{{2}} + {\frac{2}{135}}\,{u_{{6}}}^{4}{u_{{5}}}^{6}u_{{2}} -  {\frac {1}{3024}}\,{u_{{6}}}^{8}{u_{{5}}}^{4}u_{{2}}  \\
&\quad \textstyle  + {\frac {1}{158760}}\,{u_{{6}}}^{7}{u_{{5}}}^{8} -  \frac{1}{6}\,u_{{6}}{u_{{5}}}^{4}u_{{4}}u_{{1}} + {\frac {1}{1995840}}\,{u_{{6}}}^{13}{u_{{5}}}^{2}{u_{{4}}}^{2} + {\frac {1}{120}}\,{u_{{6}}}^{2}{u_{{4}}}^{7} 
+ \dots \\
\end{align*}
\begin{align*}
&\quad \textstyle \dots + {\frac {5}{18144}}\,{u_{{6}}}^{9}{u_{{2}}}^{2} + {\frac {1}{11467298304}}
\,{u_{{6}}}^{19}{u_{{5}}}^{2} + {\frac {1}{2395008}}\,{u_{{6}}}^{12}u_{{1}} + \frac{2}{3}\,u_{{6}}{u_{{5}}}^{3}u_{{4}}u_{{3}}u_{{2}} \\ 
&\quad \textstyle  + {\frac {1}{5670}}\,{u_{{6}}}^{3}{u_{{5}}}^{10} -  \frac{1}{18}\,{u_{{6}}}^{4}{u_{{5}}}^{3}u_{{3}}u_{{2}} + \frac{1}{9}\,{u_{{6}}}^{2}{u_{{5}}}^{3}{u_{{4}}}^{3}u_{{3}} -  \frac{1}{5}\,u_{{6}}{u_{{5}}}^{5}{u_{{4}}}^{2}u_{{3}} \\
&\quad \textstyle -  {\frac {43}{5987520}}\,{u_{{6}}}^{12}u_{{5}}u_{{4}}u_{{3}} + {\frac{113}{90720}}\,{u_{{6}}}^{9}{u_{{5}}}^{2}
u_{{4}}u_{{2}} -  {\frac {1}{229345966080}}\,{u_{{6}}}^{20}u_{{4}} \\
&\quad \textstyle  - {\frac {1}{90}}\,{u_{{6}}}^{5}{u_{{5}}}^{3}{u_{{4}}}^{2}u_{{3}} + {\frac{1}{540}}\,{u_{{6}}}^{5}{u_{{5}}}^{6}{u_{{4}}}^{2} + \frac{1}{2}\,{u_{{6}}}^{2}{
u_{{5}}}^{2}u_{{4}}{u_{{2}}}^{2} + {\frac {1}{17010}}\,{u_{{6}}}^{9}{u_{{5}}}^{4}{u_{{4}}}^{2} \\
&\quad \textstyle -  {\frac {1}{8064}}\,{u_{{6}}}^{8}{u_{{4}}}^{5} -  {
\frac {7}{360}}\,{u_{{5}}}^{4}{u_{{4}}}^{5}  -  {\frac {1}{3991680}}\,{u_{{6}}}^{11}{u_{{4}}}^{4} 
-  {\frac {1}{75442752}}\,{u_{{6}}}^{15}{u_{{5}}}^{4} \\
&\quad \textstyle -  {\frac {1}{199584}}\,{u_{{6}}}^{11}{u_{{3}}}^{2} + {\frac {1}{28740096}}\,{u_{{6}}}^{16}u_{{2}} + {\frac {1}{720}}\,{u_{{6}}}^{5}{u_{{4}}}^{6} -  {\frac {1}{167650560}}\,{u_{{6}}}^{14}{u_{{4}}}^{3}  \\
&\quad \textstyle - \frac{1}{12}\,{u_{{6}}}^{4}{u_{{5}}}^{2}u_{{4}}{u_{{3}}}^{2} 
-  {\frac {1}{1008}}\,{u_{{6}}}^{8}u_{{5}}u_{{3}}u_{{2}} 
+ \frac{1}{2}\,{u_{{6}}}^{2}{u_{{5}}}^{2}{u_{{4}}}^{2}u_{{1}} + {\frac{11}{90720}}\,{u_{{6}}}^{9}u_{{4}}u_{{1}} \\
&\quad \textstyle  + {\frac {11}{360}}\,{u_{{6}}}^{6}u_{{4}}{u_{{2}}}^{2} -  {u_{{6}}}^{2}u_{{5}}{u_{{4}}}^{2}u_{{3}}u_{{2}} -  u_{{6}}u_{{4}}{u_{{3}}}^{2}u_{{2}} -  u_{{6}}{u_{{5}}}^{2}{u_{{4}}}^{2}{u_{{3}}}^{2} \\
&\quad \textstyle + \frac{1}{9}\,{u_{{6}}}^{3}u_{{5}}{u_{{4}}}^{4}u_{{3}} + {\frac {11}{15120}}\,{u_{{6}}}^{8}{u_{{5}}}^{3}u_{{4}}u_{{3}} + {\frac {1}{60}}\,{u_{{6}}}^{5}{u_{{5}}}^{2}u_{{4}}u_{{1}} -  {\frac {7}{90}}\,u_{{6}}{u_{{5}}}^{6}u_{{4}}u_{{2}} \\
&\quad \textstyle -  {\frac {7}{360}}\,{u_{{6}}}^{5}{u_{{5}}}^{4}u_{{4}}u_{{2}} -  {\frac{1}{360}}\,{u_{{6}}}^{6}{u_{{5}}}^{2}{u_{{4}}}^{2}u_{{2}} -  {\frac {7}{540}}
\,{u_{{6}}}^{6}u_{{5}}{u_{{4}}}^{3}u_{{3}} - \frac{1}{30}\,{u_{{5}}}^{6}u_{{1}}  \\
&\quad \textstyle -  {\frac {1}{60}}\,{u_{{6}}}^{4}{u_{{5}}}^{5}u_
{{4}}u_{{3}} -  \frac{1}{12}\,{u_{{6}}}^{2}{u_{{5}}}^{4}{u_{{4}}}^{2}u_{{2}} + {
\frac {1}{45360}}\,{u_{{6}}}^{9}u_{{5}}{u_{{4}}}^{2}u_{{3}} + \frac{1}{3}\,{u_{{4}}}^{4}u_{{1}}\\
&\quad \textstyle + {\frac {1}{10595783632896}}\,{u_{{6}}}^{23} - \frac{1}{6}\,{u_{{6}}}^{2}{u_{{4}}}^{3}{u_{{3}}}^{2} 
+ \frac{2}{9}\,{u_{{6}}}^{3}{u_{{5}}}^{2}{u_{{4}}}^{3}u_{{2}} -  {\frac {1}{641520}}\,{u_{{6}}}^{11}{u_{{5}}}^{6}   \big]
\end{align*}
\begin{align*}
&\quad + \lambda_4^2 \cdot \big[ \textstyle   -{\frac {1}{180}}\,{u_{{6}}}^{6}{u_{{4}}}^{2}u_{{1}} - {\frac{1}{72}}\,{u_{{6}}}^{4}{u_{{4}}}^{4}u_{{2}} - {\frac {13}{2041200}}\,{u_{{6}}}^{12}{u_{{5}}}^{4}u_{{4}}  \\
&\quad \textstyle + {\frac{1}{2993760}}\,{u_{{6}}}^{13}u_{{4}}u_{{2}}+{\frac{4}{315}}\,{u_{{5}}}^{7}u_{{4}}u_{{3}} + {\frac{1}{22680}}\,{u_{{6}}}^{9}u_{{4}}u_{{1}}+\frac{1}{9}\,{u_{{6}}}^{3}{u_{{4}}}^{3}u_{{1}} \\
&\quad \textstyle + {\frac{4}{15}}\,{u_{{5}}}^{5}u_{{3}}u_{{2}} +{\frac{1}{467775}}\,{u_{{6}}}^{12}{u_{{5}}}^{2}u_{{2}} - \frac{1}{3}\,{u_{{6}}}^{3}{u_{{4}}}^{2}{u_{{2}}}^{2} - {\frac {17}{315}}\,u_{{6}}{u_{{5}}}^{8}{u_{{4}}}^{2} \\
&\quad \textstyle + {\frac{5}{756}}\,{u_{{6}}}^{4}{u_{{5}}}^{8}u_{{4}} - {\frac{1}{14175}}\,{u_{{6}}}^{8}{u_{{5}}}^{6}u_{{4}} +\frac{1}{15}\,u_{{6}}{u_{{5}}}^{2}{u_{{4}}}^{6}+{\frac{13}{120}}\,{u_{{6}}}^{4}{u_{{5}}}^{2}{u_{{4}}}^{5} \\
&\quad \textstyle  + {\frac{7}{108}}\,{u_{{6}}}^{3}{u_{{5}}}^{4}{u_{{4}}}^{4} - {\frac{1}{280}}\,{u_{{6}}}^{7}{u_{{5}}}^{2}{u_{{4}}}^{4} - \frac{1}{30}\,u_{{5}}{u_{{4}}}^{5}u_{{3}}-{\frac{1}{162}}\,{u_{{6}}}^{6}{u_{{5}}}^{4}{u_{{4}}}^{3} \\
&\quad \textstyle  + {\frac{101}{226800}}\,{u_{{6}}}^{10}{u_{{5}}}^{2}{u_{{4}}}^{3} -{\frac {4}{467775}}\,{u_{{6}}}^{11}{u_{{5}}}^{3}u_{{3}} - {\frac{8}{135}}\,{u_{{6}}}^{2}{u_{{5}}}^{6}{u_{{4}}}^{3} \\
&\quad \textstyle  - {\frac{4}{4725}}\,{u_{{6}}}^{7}{u_{{5}}}^{5}u_{{3}}  +\frac{1}{3}\,{u_{{4}}}^{3}{u_{{2}}}^{2}+{\frac {17}{838252800}}\,{u_{{6}}}^{17}{u_{{4}}}^{2}+{\frac {1}{105}}\,{u_{{5}}}^{8}u_{{2}} \\
&\quad \textstyle + \frac{1}{30}\,{u_{{6}}}^{5}{u_{{4}}}^{2}{u_{{3}}}^{2} -{\frac {1}{2520}}\,{u_{{6}}}^{8}u_{{4}}{u_{{3}}}^{2}+{\frac{17}{60}}\,u_{{6}}{u_{{4}}}^{5}u_{{2}} + {\frac{1}{16200}}\,{u_{{6}}}^{10}{u_{{4}}}^{2}u_{{2}}\\
&\quad \textstyle  +{\frac{23}{3780}}\,{u_{{6}}}^{7}{u_{{4}}}^{3}u_{{2}} - {\frac{4}{945}}\,{u_{{6}}}^{3}{u_{{5}}}^{7}u_{{3}} +\frac{1}{6}\,{u_{{5}}}^{2}{u_{{4}}}^{4}u_{{2}}+{\frac{1}{135}}\,{u_{{6}}}^{4}{u_{{5}}}^{6}u_{{2}}  \\
&\quad \textstyle -\frac{1}{3}\,{u_{{6}}}^{2}{u_{{4}}}^{3}{u_{{3}}}^{2} - {\frac {43}{396900}}\,{u_{{6}}}^{7}{u_{{5}}}^{8} -{\frac{1}{199584}}\,{u_{{6}}}^{13}{u_{{5}}}^{2}{u_{{4}}}^{2} + \frac{1}{24}\,{u_{{6}}}^{2}{u_{{4}}}^{7} \\
&\quad \textstyle - {\frac{2}{45}}\,{u_{{5}}}^{6}u_{{1}} - {\frac {1}{22680}}\,{u_{{6}}}^{9}{u_{{2}}}^{2} - {\frac{1}{28668245760}}\,{u_{{6}}}^{19}{u_{{5}}}^{2}  -{\frac{1}{14968800}}\,{u_{{6}}}^{12}u_{{1}} \\
&\quad \textstyle +{\frac{41}{14175}}\,{u_{{6}}}^{3}{u_{{5}}}^{10} - \frac{4}{9}\,{u_{{6}}}^{2}{u_{{5}}}^{3}{u_{{4}}}^{3}u_{{3}}-\frac{2}{15}\,u_{{6}}{u_{{5}}}^{5}{u_{{4}}}^{2}u_{{3}}  - {\frac {31}{7484400}}\,{u_{{6}}}^{12}u_{{5}}u_{{4}}u_{{3}} \\ 
&\quad \textstyle +{\frac{2}{2835}}\,{u_{{6}}}^{9}{u_{{5}}}^{2}u_{{4}}u_{{2}} + {\frac{1}{573364915200}}\,{u_{{6}}}^{20}u_{{4}} - {\frac {4}{45}}\,{u_{{6}}}^{5}{u_{{5}}}^{3}{u_{{4}}}^{2}u_{{3}} \\
&\quad \textstyle - {\frac {1}{1350}}\,{u_{{6}}}^{5}{u_{{5}}}^{6}{u_{{4}}}^{2} +{\frac {29}{68040}}\,{u_{{6}}}^{9}{u_{{5}}}^{4}{u_{{4}}}^{2}-{\frac {247}{100800}}
\,{u_{{6}}}^{8}{u_{{4}}}^{5}+\frac{1}{36}\,{u_{{5}}}^{4}{u_{{4}}}^{5} \\ 
&\quad \textstyle - {\frac{593}{9979200}}\,{u_{{6}}}^{11}{u_{{4}}}^{4} + {\frac {1}{67359600}}\,{u_{{6}}}^{15}{u_{{5}}}^{4} +{\frac {1}{1247400}}\,{u_{{6}}}^{11}{u_{{3}}}^{2}+{\frac{1}{2514758400}}\,{u_{{6}}}^{16}u_{{2}} \\
&\quad \textstyle + {\frac {13}{900}}\,{u_{{6}}}^{5}{u_{{4}}}^{6} +{\frac{19}{41912640}}\,{u_{{6}}}^{14}{u_{{4}}}^{3} - \frac{1}{12}\,{u_{{4}}}^{4}u_{{1}} \\
&\quad \textstyle +{\frac {1}{90}}\,{u_{{6}}}^{6}u_{{4}}{u_{{2}}}^{2} - 
{\frac{5}{18}}\,{u_{{6}}}^{3}u_{{5}}{u_{{4}}}^{4}u_{{3}} - {\frac{2}{945}}\,{u_{{6}}}^{8}{u_{{5}}}^{3}u_{{4}}u_{{3}} + {\frac{2}{45}}\,u_{{6}}{u_{{5}}}^{6}u_{{4}}u_{{2}} \\
&\quad \textstyle +{\frac {2}{45}}\,{u_{{6}}}^{5}{u_{{5}}}^{4}u_{{4}}u_{{2}}+{\frac{2}{45}}\,{u_{{6}}}^{6}{u_{{5}}}^{2}{u_{{4}}}^{2}u_{{2}} - {\frac {7}{270}}\,{u_{{6}}}^{6}u_{{5}}{u_{{4}}}^{3}u_{{3}} + 
\frac{4}{9}\,{u_{{6}}}^{3}{u_{{5}}}^{2}{u_{{4}}}^{3}u_{{2}} \\
&\quad \textstyle -{\frac {2}{45}}\,{u_{{6}}}^{4}{u_{{5}}}^{5}u_{{4}}u_{{3}} + \frac{1}{3}\,{u_{{6}}}^{2}{u_{{5}}}^{4}{u_{{4}}}^{2}u_{{2}} - {\frac {17}{22680}}\,{u_{{
6}}}^{9}u_{{5}}{u_{{4}}}^{2}u_{{3}}  \\
&\quad \textstyle - {\frac {1}{561330}}\,{u_{{6}}}^{11
}{u_{{5}}}^{6}-{\frac {1}{26489459082240}}\,{u_{{6}}}^{23} + {\frac{1}{1890}}\,{u_{{6}}}^{8}{u_{{5}}}^{4}u_{{2}} \\
&\quad \textstyle + {\frac {1}{39916800}}\,{u_{{6}}}^{16}{u_{{5}}}^{2}u_{{4}} 
- {\frac{1}{157172400}}\,{u_{{6}}}^{15}u_{{5}}u_{{3}} - {\frac{41}{4725}}\,{u_{{5}}}^{10}u_{{4}}  \big].
\end{align*}

\newpage

\section{The 4-index Q-functions} \label{APPQ}

This appendix contains a list of PDEs, expressing the 4-index $Q$-functions that were not elements of the basis for $\Gamma \big( J, \mathcal{O}(2 \Theta^{[5]} ) \big)$, as a linear combination of basis elements.   This appendix contains all the equations down to weight $-22$.  The full set can be accessed at \cite{myweb}.

The PDEs are ordered in decreasing weight (as indicated by the number in brackets).
\begin{align*}
\textbf{(-4)} \quad Q_{6666} &= - 3\wp_{55} + 4\wp_{46}, \\
\textbf{(-5)} \quad Q_{5666} &= - 2\wp_{45}, \\
\textbf{(-6)} \quad Q_{4666} &= 6\lambda_{4}\wp_{66} - 2\wp_{44} - \textstyle \frac{3}{2}Q_{5566}, \\
\textbf{(-7)} \quad Q_{4566} &= 2\lambda_{4}\wp_{56} + 2\wp_{36}, \\
\textbf{(-7)} \quad Q_{5556} &= 4\lambda_{4}\wp_{56} + 4\wp_{36}, \\
\textbf{(-8)} \quad Q_{4466} &= 4\lambda_4\wp_{46} + \lambda_{4}\wp_{55} + \wp_{35} + 6\wp_{26} - Q_{4556} + 4\lambda_3, \\
\textbf{(-8)} \quad Q_{5555} &= 16\lambda_{4}\wp_{55} + 4\wp_{35} + 24\wp_{26} - 6Q_{4556} + 16\lambda_3, \\
\textbf{(-9)} \quad Q_{4456} &= \textstyle \frac{8}{3}\lambda_{4}\wp_{45} + 2\wp_{25} - \frac{4}{3}\wp_{34} -  \frac{1}{6}Q_{4555}, \\
\textbf{(-9)} \quad Q_{3666} &= \textstyle 2\lambda_{4}\wp_{45} - \frac{1}{2}Q_{4555}, \\
\textbf{(-10)} \quad Q_{2666} &= \textstyle - \frac{1}{2}Q_{4455} - \frac{1}{2}Q_{3566} - \frac{1}{2}\lambda_{4}Q_{5566} + 4\wp_{66}\lambda_3, \\
\textbf{(-10)} \quad Q_{4446} &= \textstyle 6\lambda_{4}\wp_{44} - 2\wp_{24} - 2\lambda_{4}^2\wp_{66} \\
&\quad \textstyle - Q_{4455} + \frac{1}{2}\lambda_{4}Q_{5566} - \frac{5}{2}Q_{3566} + 4\wp_{66}\lambda_{3}, \\
\textbf{(-11)} \quad Q_{3466} &= \textstyle 4\wp_{36}\lambda_4 - \frac{1}{2}Q_{3566}, \\
\textbf{(-11)} \quad Q_{4445} &= \textstyle 6\wp_{36}\lambda_4 - 2\wp_{56}\lambda_4^2 + 8\wp_{56}\lambda_3 -  3Q_{2566} - \frac{3}{2}Q_{3556}, \\
\textbf{(-12)} \quad Q_{2466} &= \textstyle 8\wp_{16} - \wp_{33} + 2\lambda_2 - \frac{1}{2}Q_{2456} - Q_{3456} + 4\lambda_4\wp_{26} + 2\lambda_4\wp_{35}, \\
\textbf{(-12)} \quad Q_{3555} &= 24\wp_{16} - 8\wp_{33} + 16\lambda_4\wp_{35} - 6Q_{3456}, \\
\textbf{(-12)} \quad Q_{4444} &=  -12\wp_{16} + 9\wp_{33} + 6\lambda_2 + 12Q_{3456} + 12\wp_{55}\lambda_3 
+ 4\wp_{46}\lambda_4^2 \\ &\quad - 16\wp_{46}\lambda_3 - 3\wp_{55}\lambda_4^2 + 12\lambda_4\wp_{26} - 18\lambda_4\wp_{35} - 6Q_{2556}, \\
\textbf{(-13)} \quad Q_{2555} &= \textstyle 2\wp_{15} - 8\wp_{23} + 16\wp_{25}\lambda_4 - 6Q_{2456}, \\
\textbf{(-13)} \quad Q_{3446} &= \textstyle 6\wp_{15} - 4\wp_{23} + 4\wp_{25}\lambda_4 + \frac{8}{3}\wp_{34}\lambda_4 \\
&\quad \textstyle + \frac{2}{3}\wp_{45}\lambda_4^2 - \frac{1}{6}Q_{4555}\lambda_4 - 2Q_{2456}, \\
\textbf{(-13)} \quad Q_{3455} &= \textstyle -4\wp_{25}\lambda_4 + \frac{8}{3}\wp_{34}\lambda_4  - \frac{4}{3}\wp_{45}\lambda_4^2 + \frac{1}{3}Q_{4555}\lambda_4 + 2Q_{2456}, \\
\textbf{(-14)} \quad Q_{1666} &= \textstyle - \frac{1}{2}\lambda_{4}Q_{3566} + Q_{3366} 
- \frac{1}{2}Q_{3445}, \\
\textbf{(-14)} \quad Q_{2446} &= \textstyle \frac{8}{3}\wp_{14} - 2\wp_{22} + \frac{8}{3}\wp_{24}\lambda_{4} 
- \frac{2}{3}\wp_{66}\lambda_{4}\lambda_{3} + \frac{5}{6}Q_{3366} + Q_{5566}\lambda_{3} \\
&\quad \textstyle - \frac{1}{6}Q_{5566}\lambda_{4}^{2} 
+ 2\wp_{44}\lambda_{3} - \frac{1}{6}\lambda_{4}Q_{4455} - \frac{4}{3}\lambda_{4}Q_{3566}  - \frac{7}{6}Q_{3445}, \\
\textbf{(-14)} \quad Q_{2455} &= \textstyle \frac{8}{3}\wp_{14} - 2\wp_{22} + \frac{8}{3}\wp_{24}\lambda_{4} 
- \frac{2}{3}\wp_{66}\lambda_{4}\lambda_{3} - \frac{5}{3}Q_{3366} - Q_{5566}\lambda_{3} \\
&\quad \textstyle + \frac{1}{3}Q_{5566}\lambda_{4}^{2} 
+ 6\wp_{66}\lambda_{2} - 2\wp_{44}\lambda_{3} + \frac{1}{3}\lambda_{4}Q_{4455} 
+ \frac{2}{3}\lambda_{4}Q_{3566} \\
&\quad \textstyle + \frac{1}{3}Q_{3445},
\end{align*}

\newpage

\begin{align*}
\textbf{(-15)} \quad Q_{1566} &= \textstyle - \frac{1}{4}Q_{2445} - \frac{1}{4}Q_{2455}\lambda_4 + \frac{3}{4}\wp_{2366} - \frac{3}{2}\wp_{23}\wp_{66} - 3\wp_{26}\wp_{36} \\
&\quad \textstyle - \frac{1}{2}\wp_{36}\lambda_3 + \frac{3}{2}\wp_{56}\lambda_2, \\
\textbf{(-15)} \quad Q_{3356} &= \textstyle - \frac{1}{2}Q_{2445} + 3\wp_{36}\lambda_3 + 3\wp_{56}\lambda_2 
- \frac{1}{2}Q_{2566}\lambda_4 - \frac{1}{2}Q_{2366}, \\
\textbf{(-15)} \quad Q_{3444} &= \textstyle - \frac{3}{4}Q_{2445} + 4\wp_{36}\lambda_4^2 - \frac{3}{2}Q_{3556}\lambda_4 
+ \frac{13}{2}\wp_{36}\lambda_3 + \frac{9}{2}\wp_{56}\lambda_2 \\
&\quad \textstyle - \frac{3}{4}Q_{2566}\lambda_4 - \frac{15}{4}Q_{2366}, \\
\textbf{(-16)} \quad Q_{3346} &= -2Q_{1466} + 4\lambda_1 - Q_{1556} - Q_{2266} + 16\wp_{16}\lambda_4 
+ 3\wp_{35}\lambda_3 - \wp_{55}\lambda_2 \\
&\quad \textstyle + 2\wp_{26}\lambda_3 - 2Q_{2356}, \\
\textbf{(-16)} \quad Q_{3355} &= 8Q_{1466} - 8\lambda_1 + 6Q_{1556} - 32\wp_{16}\lambda_4 + 4\wp_{35}\lambda_3 \\
&\quad \textstyle - 2Q_{2356}, \\
\textbf{(-16)} \quad Q_{2444} &= \textstyle + 2\lambda_{4}\lambda_{2} + 10\lambda_{1} + 36\wp_{16}\lambda_{4} 
+ \frac{3}{2}\wp_{35}\lambda_{3} + \frac{15}{2}\wp_{55}\lambda_{2} - Q_{3456}\lambda_{4} \\
&\quad \textstyle - \frac{3}{2}Q_{2556}\lambda_{4} - 7\wp_{26}\lambda_{3} - 6\wp_{46}\lambda_{2} + 3\wp_{33}\lambda_{4} 
+ 6\lambda_{4}\wp_{3456} \\
&\quad \textstyle - 12\lambda_{4}\wp_{34}\wp_{56} - 12\lambda_{4}\wp_{35}\wp_{46} - 12\lambda_{4}\wp_{36}\wp_{45} 
+ 4\wp_{26}\lambda_{4}^{2} - 6\wp_{35}\lambda_{4}^{2} \\
&\quad \textstyle - 9Q_{1466} - 9Q_{1556} + \frac{3}{2}Q_{2266} + 3Q_{2356}, \\
\textbf{(-17)} \quad Q_{1456} &= \textstyle - 2\wp_{13} - \wp_{45}\lambda_2 + 2\wp_{25}\wp_{26} + \wp_{22}\wp_{56} 
+ 4\wp_{15}\lambda_4 - Q_{2346} \\
&\quad \textstyle - \frac{1}{12}Q_{4555}\lambda_3 - \frac{1}{2}\wp_{2256} + \wp_{25}\lambda_3
+ \frac{1}{3}\wp_{45}\lambda_4\lambda_3 + \frac{4}{3}\wp_{34}\lambda_3, \\
\textbf{(-17)} \quad Q_{1555} &= \textstyle 4\wp_{13} + 6\wp_{45}\lambda_2 + 3Q_{2256} - 8\wp_{15}\lambda_4
+ 6Q_{2346} + \frac{1}{2}Q_{4555}\lambda_3 \\
&\quad \textstyle - 6\wp_{25}\lambda_3 - 2\wp_{45}\lambda_4\lambda_3 - 8\wp_{34}\lambda_3, \\
\textbf{(-17)} \quad Q_{2355} &= \textstyle - 4\wp_{13} - 4\wp_{45}\lambda_2 - 2Q_{2256} + 8\wp_{15}\lambda_4 
- 2Q_{2346} \\ &\quad \textstyle + 4\wp_{25}\lambda_3 + 4\wp_{34}\lambda_3, \\
\textbf{(-17)} \quad Q_{3345} &= \textstyle - 2\wp_{25}\lambda_3 + \frac{4}{3}\wp_{34}\lambda_3 
- \frac{2}{3}\wp_{45}\lambda_4\lambda_3 + 2\wp_{45}\lambda_2  + Q_{2256} + \frac{1}{6}Q_{4555}\lambda_3, \\
\textbf{(-18)} \quad Q_{1446} &= \textstyle - 6\wp_{12} + 6\wp_{66}\lambda_{1} + \frac{1}{2}Q_{3566}\lambda_{4}^2 - Q_{1455}
+ 8\wp_{14}\lambda_{4} - Q_{3566}\lambda_{3} \\
&\quad \textstyle - Q_{2345} - \frac{1}{2}Q_{3344} - \frac{1}{2}Q_{3366}\lambda_{4} + \frac{1}{2}Q_{3445}\lambda_{4}, \\
\textbf{(-18)} \quad Q_{2246} &= \textstyle 2\wp_{44}\lambda_{2} - 4\wp_{66}\lambda_{1} 
- 2\wp_{66}\lambda_{4}\lambda_{2} + \frac{1}{2}Q_{3566}\lambda_{3} + Q_{3366}\lambda_{4} \\
&\quad \textstyle - Q_{3445}\lambda_{4} - Q_{3566}\lambda_{4}^2 
+ 2\wp_{24}\lambda_{3} + \frac{3}{2}Q_{5566}\lambda_{2} + Q_{1455}, \\
\textbf{(-18)} \quad Q_{2255} &= \textstyle - 8\wp_{12} - 8\wp_{44}\lambda_{2} + 16\wp_{66}\lambda_{1} 
+ \frac{32}{3}\wp_{14}\lambda_{4} + 8\wp_{66}\lambda_{4}\lambda_{2} - 2\wp_{2345} \\
&\quad \textstyle - \frac{2}{3}Q_{3566}\lambda_{3} - \frac{8}{3}Q_{3366}\lambda_{4} + \frac{4}{3}Q_{3445}\lambda_{4} 
+ \frac{4}{3}Q_{3566}\lambda_{4}^2 + \frac{1}{3}Q_{4455}\lambda_{3} \\
&\quad \textstyle + \frac{1}{3}\lambda_{4}\lambda_{3}\wp_{5566} - \frac{2}{3}\lambda_{4}\lambda_{3}\wp_{55}\wp_{66} 
- \frac{4}{3}\lambda_{4}\lambda_{3}\wp_{56}^2 - \frac{8}{3}\wp_{66}\lambda_{3}^2 + \frac{8}{3}\wp_{24}\lambda_{3} \\
&\quad \textstyle + 4\wp_{23}\wp_{45} + 4\wp_{24}\wp_{35} + 4\wp_{25}\wp_{34} - 3Q_{5566}\lambda_{2} - 2Q_{1455},  \\
\textbf{(-19)} \quad Q_{1366} &= \textstyle - \frac{3}{10}Q_{2245} - \frac{2}{5}Q_{2344} 
- \frac{3}{10}Q_{2566}\lambda_{3} - \frac{1}{4}Q_{2366}\lambda_{4} + \frac{3}{20}Q_{2455}\lambda_{4}  \\
&\quad \textstyle + \frac{3}{20}Q_{2566}\lambda_{4}^2 + \frac{12}{5}\wp_{56}\lambda_{1} - \frac{3}{10}\wp_{56}\lambda_{4}\lambda_{2} + \frac{7}{10}\wp_{36}\lambda_{4}\lambda_{3} \\
&\quad \textstyle + \frac{9}{5}\wp_{36}\lambda_{2} - \frac{1}{4}Q_{3556}\lambda_{3}, \\
\textbf{(-19)} \quad Q_{3336} &= \textstyle - \frac{3}{5}Q_{2245} + \frac{6}{5}Q_{2344} 
- \frac{3}{5}Q_{2566}\lambda_{3} + \frac{3}{2}Q_{2366}\lambda_{4} + \frac{3}{10}Q_{2455}\lambda_{4} \\
&\quad \textstyle + \frac{3}{10}Q_{2566}\lambda_{4}^2 + \frac{38}{5}\wp_{36}\lambda_{2} - \frac{16}{5}\wp_{56}\lambda_{1} - \frac{3}{5}\wp_{56}\lambda_{4}\lambda_{2} \\
&\quad \textstyle - \frac{33}{5}\wp_{36}\lambda_{4}\lambda_{3} + \frac{3}{2}Q_{3556}\lambda_{3}, 
\end{align*}

\newpage

\begin{align*}
\textbf{(-19)} \quad Q_{1445} &= \textstyle - \frac{1}{2}Q_{2245} - Q_{2344} - \frac{1}{2}Q_{2566}\lambda_{3} 
- \frac{3}{2}Q_{2366}\lambda_{4}  + 8\wp_{56}\lambda_{1} + 3\wp_{36}\lambda_{2} \\
&\quad \textstyle + \frac{1}{2}Q_{2566}\lambda_{4}^2   + \frac{1}{2}Q_{2455}\lambda_{4}
- 2\wp_{56}\lambda_{4}\lambda_{2} + 3\wp_{36}\lambda_{4}\lambda_{3} - \frac{3}{4}Q_{3556}\lambda_{3}, \\
\textbf{(-20)} \quad Q_{2244} &= \textstyle + 2\lambda_{3}\lambda_{2} + 28\lambda_{4}\lambda_{1} + 2\wp_{55}\lambda_{4}\lambda_{2} 
- 6\wp_{35}\lambda_{4}\lambda_{3} + 3\wp_{33}\lambda_{3}  \\
&\quad \textstyle - 24\wp_{46}\lambda_{1} + 112\wp_{16}{\lambda_{4}}^{2} - 44\wp_{16}\lambda_{3} + 28\wp_{55}\lambda_{1} 
+ 9\wp_{35}\lambda_{2}  \\
&\quad \textstyle + 20Q_{1356} - 6Q_{1444} + 10Q_{1266} - 22Q_{1466}\lambda_{4} + Q_{2266}\lambda_{4}  \\
&\quad \textstyle + 2Q_{2356}\lambda_{4} + 2Q_{3456}\lambda_{3} - 22Q_{1556}\lambda_{4} - Q_{2556}\lambda_{3} - 6\wp_{26}\lambda_{2}, \\
\textbf{(-20)} \quad Q_{2336} &= \textstyle 8\wp_{16}\lambda_{3} - 2\wp_{55}\lambda_{1} + 2\wp_{35}\lambda_{2} - 2Q_{1356} - 2Q_{1266}, \\
\textbf{(-20)} \quad Q_{3335} &= \textstyle - 32\lambda_{4}\lambda_{1} + 24\wp_{46}\lambda_{1} - 128\wp_{16}{\lambda_{4}}^{2} 
+ 56\wp_{16}\lambda_{3} - 32\wp_{55}\lambda_{1}  \\
&\quad \textstyle - 30Q_{1356} + 8Q_{1444} - 12Q_{1266} + 24Q_{1466}\lambda_{4} + 24Q_{1556}\lambda_{4} \\
&\quad \textstyle + 4\wp_{35}\lambda_{2}, \\
\textbf{(-21)} \quad Q_{1256} &= \textstyle - \frac{1}{2}Q_{2236} - Q_{1346} - \frac{1}{12}Q_{4555}\lambda_{2} 
+ \frac{1}{3}\wp_{45}\lambda_{4}\lambda_{2} + 3\wp_{15}\lambda_{3} \\
&\quad \textstyle + \frac{4}{3}\wp_{34}\lambda_{2} 
- 2\wp_{45}\lambda_{1}, \\
\textbf{(-21)} \quad Q_{1355} &= \textstyle \frac{1}{2}Q_{2236} + Q_{1346} + \frac{1}{4}Q_{4555}\lambda_{2} - \frac{1}{2}Q_{2335} - \wp_{45}\lambda_{4}\lambda_{2} + \wp_{15}\lambda_{3}, \\
\textbf{(-21)} \quad Q_{3334} &= \textstyle - \frac{3}{2}Q_{2236} - 3Q_{1346} + \frac{1}{4}Q_{4555}\lambda_{2} 
+ \frac{3}{2}Q_{2335} - \wp_{45}\lambda_{4}\lambda_{2} + 9\wp_{15}\lambda_{3} \\
&\quad \textstyle + 12\wp_{34}\lambda_{2} - 8\wp_{45}\lambda_{1}, \\
\textbf{(-22)} \quad Q_{1345} &= \textstyle  5\wp_{{66}}\lambda_{{0}}
- \wp_{11} + \frac{16}{15}\wp_{14}\lambda_4^2 
+ \frac{11}{5}\wp_{{66}}\lambda_{{4}}\lambda_{{1}} - {\frac{7}{30}}Q_{3566}\lambda_{4}\lambda_{3}  \\
&\quad \textstyle   - \frac{1}{5} Q_{3344}\lambda_{{4}}
- \frac{1}{2} Q_{1255} + \frac{10}{3}\wp_{{14}}\lambda_{{3}} - \wp_{{44}}\lambda_{{1}} -  Q_{1246}  \\
&\quad \textstyle + \frac{2}{15} Q_{3566}{\lambda_{{4}}}^{3}  
- \frac{1}{15} Q_{3366}{\lambda_{{4}}}^{2} - \frac{1}{3} Q_{3366}\lambda_{{3}} 
+ \frac{1}{6} Q_{3445}\lambda_{{3}} \\
&\quad \textstyle - \frac{2}{5} Q_{2345}\lambda_{{4}} - \frac{4}{5}\wp_{12}\lambda_{4}
+ \frac{2}{15} Q_{3445}{\lambda_{{4}}}^{2}, \\
\textbf{(-22)} \quad Q_{2226} &= \textstyle  - 6\wp_{{11}} + \frac{7}{2} Q_{3566}\lambda_{{2}} 
+ \frac{1}{2} Q_{4455}\lambda_{{2}} + 6 Q_{5566}\lambda_{{1}} 
+ {\frac {32}{5}}\wp_{{14}}{\lambda_{{4}}}^{2} \\
&\quad \textstyle - {\frac {14}{5}}\wp_{{66}}\lambda_{{4}}\lambda_{{1}} 
- {\frac {39}{10}} Q_{3566}\lambda_{{4}}\lambda_{{3}} - {\frac {24}{5}}\wp_{{12}}\lambda_{{4}} 
+ 4\wp_{{24}}\lambda_{{2}} - 10\wp_{{66}}\lambda_{{0}} \\
&\quad \textstyle + 12\wp_{{14}}\lambda_{{3}} 
+ 6\wp_{{44}}\lambda_{{1}} - 2\wp_{{66}}\lambda_{{3}}\lambda_{{2}} - 6 Q_{1246} 
+ \frac{4}{5} Q_{3445}{\lambda_{{4}}}^{2} \\
&\quad \textstyle - \frac{6}{5} Q_{3344}\lambda_{{4}} 
- \frac{2}{5} Q_{3366}{\lambda_{{4}}}^{2} + \frac{3}{2} Q_{3366}\lambda_{{3}} 
- \frac{3}{2} Q_{3445}\lambda_{{3}} - {\frac {12}{5}} Q_{2345}\lambda_{{4}} \\
&\quad \textstyle + \frac{4}{5} Q_{3566}{\lambda_{{4}}}^{3} + \frac{1}{2} Q_{5566}\lambda_{{4}}\lambda_{{2}}, \\
\textbf{(-22)} \quad Q_{2235} &= \textstyle  - 2\wp_{{11}} -  Q_{3566}\lambda_{{2}} 
- 4 Q_{5566}\lambda_{{1}} - {\frac {32}{15}}\wp_{{14}}{\lambda_{{4}}}^{2} 
+ {\frac {38}{5}}\wp_{{66}}\lambda_{{4}}\lambda_{{1}} \\
&\quad \textstyle + {\frac {22}{15}} Q_{3566}\lambda_{{4}}\lambda_{{3}} 
+ \frac{8}{5}\wp_{{12}}\lambda_{{4}} + 10\wp_{{66}}\lambda_{{0}} -  Q_{1255} 
+ \frac{4}{3}\wp_{{14}}\lambda_{{3}} \\
&\quad \textstyle - 10\wp_{{44}}\lambda_{{1}} + 2 Q_{1246} 
- {\frac {4}{15}} Q_{3445}{\lambda_{{4}}}^{2} + \frac{2}{5} Q_{3344}\lambda_{{4}} 
+ \frac{2}{15} Q_{3366}{\lambda_{{4}}}^{2} \\
&\quad \textstyle - \frac{4}{3} Q_{3366}\lambda_{{3}} 
+ \frac{2}{3} Q_{3445}\lambda_{{3}} + \frac{4}{5} Q_{2345}\lambda_{{4}} 
- {\frac {4}{15}} Q_{3566}{\lambda_{{4}}}^{3}, \\
\textbf{(-22)} \quad Q_{2334} &= \textstyle  - 4\wp_{{11}} -  Q_{3566}\lambda_{{2}} 
+ 2 Q_{5566}\lambda_{{1}} - {\frac {32}{15}}\wp_{{14}}{\lambda_{{4}}}^{2} 
- {\frac {32}{5}}\wp_{{66}}\lambda_{{4}}\lambda_{{1}} \\
&\quad \textstyle + {\frac {7}{15}} Q_{3566}\lambda_{{4}}\lambda_{{3}} 
+ \frac{8}{5}\wp_{{12}}\lambda_{{4}} +  Q_{1255} + \frac{4}{3}\wp_{{14}}\lambda_{{3}} 
+ 4\wp_{{44}}\lambda_{{1}} \\
&\quad \textstyle - {\frac {4}{15}} Q_{3445}{\lambda_{{4}}}^{2} 
+ \frac{2}{5} Q_{3344}\lambda_{{4}} + \frac{2}{15} Q_{3366}{\lambda_{{4}}}^{2} 
+ \frac{2}{3} Q_{3366}\lambda_{{3}}  \\
&\quad \textstyle - \frac{1}{3} Q_{3445}\lambda_{{3}} + \frac{4}{5} Q_{2345}\lambda_{{4}} - {\frac {4}{15}}Q_{3566}{\lambda_{{4}}}^{3}  \\
\end{align*}

\newpage

\section{The two-term addition formula}\label{APPadd}

The Abelian functions associated with $C$ satisfy the following two-term addition formula.
\[
- \frac{\s(\bm{u}+\bm{v})\s(\bm{u}-\bm{v})}{\s(\bm{u})^2\s(\bm{v})^2} = f(\bm{u},\bm{v}) - f(\bm{v},\bm{u})
\]
where $f(\bu,\bv)$ is a finite polynomial of Abelian functions.  We write $f(\bu,\bv)$ as 
\[
f(\bu,\bv) = P_{30} + P_{26} + P_{22} + P_{18} + P_{14} + P_{10} + P_{6} + P_{2}
\]
where each $P_{k}$ contains the terms with weight $-k$ in the Abelian functions and weight $k-30$ in $\bm{\lambda}$.
\begin{align*}
&P_{30} = \textstyle \frac{1}{4}Q_{114466}(\bu) + \frac{5}{3}Q_{3566}(\bv)Q_{1356}(\bu) 
- \frac{1}{3}\wp_{14}(\bv)Q_{2356}(\bu) \\
&\quad \textstyle + \frac{1}{3}Q_{1556}(\bv)\wp_{14}(\bu) - \frac{1}{4}Q_{2236}(\bu)\wp_{25}(\bv) 
- \frac{2}{5}Q_{2345}(\bu)\wp_{33}(\bv) \\
&\quad \textstyle + \frac{1}{10}Q_{3344}(\bv)Q_{3456}(\bu) - \frac{1}{2}\wp_{25}(\bv)Q_{1346}(\bu) 
+ \frac{1}{2}Q_{4556}(\bv)\wp_{11}(\bu) \\
&\quad \textstyle + \frac{3}{5}Q_{3344}(\bu)\wp_{16}(\bv) + \frac{1}{24}Q_{1346}(\bv)Q_{4555}(\bu) 
+ \frac{1}{2}Q_{1145}(\bv)\wp_{56}(\bu) \\
&\quad \textstyle + \frac{1}{12}Q_{2236}(\bv)\wp_{34}(\bu) + \frac{1}{2}Q_{1444}(\bv)Q_{3566}(\bu) 
- \wp_{46}(\bv)Q_{1146}(\bu) \\
&\quad \textstyle - \frac{1}{2}Q_{1556}(\bv)\wp_{22}(\bu) + \frac{1}{2}Q_{1155}(\bu)\wp_{46}(\bv) 
+ \frac{5}{6}Q_{3566}(\bv)Q_{1266}(\bu)  \\
&\quad \textstyle + \frac{7}{5}\wp_{12}(\bv)Q_{3456}(\bu) + \wp_{26}(\bv)Q_{1246}(\bu) + Q_{2346}(\bv)\wp_{15}(\bu) 
- \frac{1}{5}\wp_{12}(\bv)Q_{2556}(\bu)  \\
&\quad \textstyle + \frac{1}{10}Q_{2345}(\bu)Q_{2556}(\bv) - \frac{2}{3}Q_{1356}(\bv)Q_{4455}(\bu) 
- \frac{1}{20}Q_{2245}(\bv)Q_{2566}(\bu) \\
&\quad \textstyle + \wp_{36}(\bv)Q_{1245}(\bu)  - \frac{1}{5}\wp_{33}(\bv)Q_{3344}(\bu) 
+ Q_{1246}(\bv)\wp_{35}(\bu) - \frac{2}{3}\wp_{24}(\bv)Q_{1444}(\bu) \\
&\quad \textstyle + \frac{1}{3}Q_{1466}(\bu)Q_{3445}(\bv) - \frac{4}{3}Q_{1466}(\bu)\wp_{14}(\bv) 
- \wp_{26}(\bv)\wp_{11}(\bu) - Q_{2346}(\bv)\wp_{23}(\bu) \\
&\quad \textstyle + \frac{1}{5}\wp_{16}(\bv)Q_{2345}(\bu) - \frac{1}{12}Q_{3445}(\bv)Q_{2266}(\bu) 
+ \frac{1}{2}Q_{5566}(\bv)Q_{1244}(\bu) \\
&\quad \textstyle + \frac{1}{6}Q_{1466}(\bv)Q_{3366}(\bu) + \frac{1}{4}8Q_{4555}(\bv)Q_{2335}(\bu) 
+ \frac{1}{3}Q_{4455}(\bv)Q_{1266}(\bu) \\
&\quad \textstyle - \frac{1}{6}Q_{4455}(\bv)Q_{1444}(\bu) + \frac{1}{5}Q_{2345}(\bv)Q_{3456}(\bu) 
+ \frac{1}{2}Q_{1156}(\bv)\wp_{45}(\bu) \\
&\quad \textstyle  + Q_{2256}(\bv)\wp_{15}(\bu) - \frac{32}{5}\wp_{12}(\bv)\wp_{16}(\bu)
- \frac{7}{3}Q_{1356}(\bv)\wp_{24}(\bu) + \frac{1}{4}Q_{2335}(\bv)\wp_{25}(\bu) \\
&\quad \textstyle - \frac{1}{2}\wp_{33}(\bv)Q_{1455}(\bu) + \frac{1}{2}\wp_{44}(\bv)Q_{1166}(\bu) 
- \wp_{35}(\bv)\wp_{11}(\bu) - Q_{1455}(\bv)\wp_{16}(\bu) \\
&\quad \textstyle - \frac{1}{48}Q_{2236}(\bv)Q_{4555}(\bu)  - \frac{1}{20}Q_{3344}(\bv)Q_{2556}(\bu) 
- \frac{1}{10}Q_{2566}(\bv)Q_{2344}(\bu) \\
&\quad \textstyle - \frac{1}{10}Q_{3556}(\bv)Q_{2245}(\bu) - \frac{1}{20}Q_{3556}(\bv)Q_{2344}(\bu) 
+ \frac{1}{8}Q_{2445}(\bv)Q_{2366}(\bu) \\
&\quad \textstyle + \frac{1}{6}Q_{1556}(\bv)Q_{3366}(\bu) - \frac{1}{2}Q_{1255}(\bu)\wp_{26}(\bv) 
- \frac{1}{2}Q_{1144}(\bv)\wp_{66}(\bu) \\
&\quad \textstyle - \frac{5}{12}Q_{2335}(\bv)\wp_{34}(\bu)  + \frac{1}{6}Q_{3366}(\bv)Q_{2266}(\bu)
+ \frac{2}{3}\wp_{24}(\bv)Q_{1266}(\bu) \\
&\quad \textstyle - \frac{1}{2}Q_{2256}(\bv)\wp_{23}(\bu) - \wp_{22}(\bu)Q_{1466}(\bv) 
+ \frac{1}{2}\wp_{55}(\bv)Q_{1146}(\bu) + \wp_{44}(\bv)Q_{1244}(\bu)  \\
&\quad \textstyle + \frac{1}{6}Q_{1346}(\bu)\wp_{34}(\bv)  - \frac{1}{3}\wp_{14}(\bu)Q_{2266}(\bv) 
+ \wp_{13}(\bv)Q_{2456}(\bu) + \frac{14}{5}\wp_{33}(\bu)\wp_{12}(\bv) \\
&\quad \textstyle - \frac{1}{3}Q_{2356}(\bv)Q_{3366}(\bu) + \frac{1}{3}Q_{3445}(\bv)Q_{1556}(\bu) 
- \frac{1}{6}Q_{3445}(\bv)Q_{2356}(\bu).
\end{align*}
\begin{align*}
&P_{26} = \textstyle \big[ \frac{1}{4}Q_{1155}(\bu) - Q_{1146}(\bu) - \frac{1}{60}Q_{2556}(\bu)Q_{3366}(\bv) 
- \frac{1}{2}Q_{1556}(\bu)Q_{4455}(\bv) \\
&\quad \textstyle + \frac{1}{24}Q_{3566}(\bu)Q_{2266}(\bv) - \frac{7}{6}Q_{5566}(\bv)Q_{1266}(\bu) 
+ \frac{1}{2}Q_{1466}(\bv)Q_{4455}(\bu) \\
&\quad \textstyle - \frac{1}{15}Q_{3445}(\bv)Q_{3456}(\bu) - \frac{7}{3}Q_{5566}(\bv)Q_{1356}(\bu) 
+ \frac{7}{6}Q_{3566}(\bu)Q_{1556}(\bv) \\
&\quad \textstyle + \frac{1}{30}Q_{2556}(\bu)Q_{3445}(\bv) - \frac{1}{2}Q_{5566}(\bu)Q_{1444}(\bv) + \frac{1}{20}Q_{3556}(\bv)Q_{2445}(\bu) \\
&\quad \textstyle - \frac{1}{30}Q_{3456}(\bv)Q_{3366}(\bu) + \frac{1}{40}Q_{2566}(\bu)Q_{2445}(\bv) 
+ \frac{7}{6}Q_{1466}(\bv)Q_{3566}(\bu) \\
&\quad \textstyle + \frac{1}{6}Q_{3566}(\bu)Q_{2356}(\bv) - \wp_{11}(\bu)\wp_{55}(\bv) 
- \frac{3}{10}Q_{3445}(\bv)\wp_{33}(\bu) \\
&\quad \textstyle + \frac{4}{15}Q_{2556}(\bu)\wp_{14}(\bv) - \frac{1}{6}\wp_{45}(\bu)Q_{1346}(\bv) 
- \frac{128}{15}\wp_{16}(\bv)\wp_{14}(\bu)  \\
&\quad \textstyle + 2\wp_{23}(\bu)\wp_{15}(\bv) - \frac{8}{15}Q_{3366}(\bv)\wp_{16}(\bu) 
- \frac{1}{24}Q_{2335}(\bu)\wp_{45}(\bv)  \\
&\quad \textstyle + \frac{2}{5}\wp_{33}(\bu)Q_{3366}(\bv) + 2Q_{1466}(\bv)\wp_{24}(\bu) + 4\wp_{16}(\bv)\wp_{22}(\bu) 
- 2Q_{1556}(\bu)\wp_{24}(\bv) \\
&\quad \textstyle + 3\wp_{44}(\bu)Q_{1266}(\bv)  + \frac{4}{3}\wp_{44}(\bv)Q_{1444}(\bu) 
+ \frac{3}{5}\wp_{36}(\bv)Q_{2344}(\bu)  \\
&\quad \textstyle - 2\wp_{35}(\bu)\wp_{12}(\bv) + \frac{1}{24}\wp_{45}(\bu)Q_{2236}(\bv) 
- \frac{8}{15}Q_{3456}(\bu)\wp_{14}(\bv)  \\
&\quad \textstyle - \wp_{66}(\bv)Q_{1244}(\bu) - 2\wp_{25}(\bu)\wp_{13}(\bv) - \frac{1}{5}\wp_{36}(\bu)Q_{2245}(\bv) \\
&\quad \textstyle + \frac{14}{15}\wp_{16}(\bv)Q_{3445}(\bu) + 6\wp_{44}(\bu)Q_{1356}(\bv) 
- \frac{12}{5}\wp_{33}(\bu)\wp_{14}(\bv) \big] \lambda_4
\end{align*}
\begin{align*}
&P_{22} = \textstyle \big[ \wp_{35}(\bv)\wp_{14}(\bu) - Q_{1466}(\bv)Q_{5566}(\bu) 
- \wp_{23}(\bv)\wp_{25}(\bu) - \frac{58}{15}\wp_{16}(\bu)Q_{3566}(\bv)  \\
&\quad \textstyle - Q_{1556}(\bu)Q_{5566}(\bv)  - \frac{1}{40}Q_{2566}(\bu)Q_{3556}(\bv) 
- \frac{1}{5}Q_{3456}(\bu)Q_{3566}(\bv) \\
&\quad \textstyle - 2Q_{1466}(\bu)\wp_{44}(\bv) - \frac{1}{2}Q_{3366}(\bv)\wp_{35}(\bu) 
+ 2\wp_{44}(\bu)Q_{1556}(\bv)  \\
&\quad \textstyle + \frac{1}{2}\wp_{15}(\bv)\wp_{25}(\bu) 
- \frac{1}{24}\wp_{23}(\bu)Q_{4555}(\bv) + \frac{2}{5}\wp_{33}(\bu)Q_{3566}(\bv) \\
&\quad \textstyle + \frac{8}{3}Q_{1266}(\bu)\wp_{66}(\bv) + \frac{16}{3}Q_{1356}(\bv)\wp_{66}(\bu) 
- \frac{4}{3}\wp_{23}(\bv)\wp_{34}(\bu) \\
&\quad \textstyle + \frac{1}{10}Q_{2556}(\bu)Q_{3566}(\bv) 
+ \frac{4}{3}Q_{1444}(\bu)\wp_{66}(\bv) + \frac{1}{4}Q_{3445}(\bv)\wp_{35}(\bu) \\
&\quad \textstyle - \frac{4}{3}\wp_{16}(\bu)Q_{4455}(\bv)  - 4\wp_{16}(\bu)\wp_{24}(\bv) 
- \frac{1}{4}Q_{2366}(\bu)\wp_{36}(\bv) - \frac{13}{6}\wp_{15}(\bu)\wp_{34}(\bv) \\
&\quad \textstyle + \frac{1}{24}\wp_{15}(\bu)Q_{4555}(\bv) \big] \lambda_3
+ \big[ - \frac{1}{15}Q_{3456}(\bu)Q_{3566}(\bv) - 4\wp_{44}(\bu)Q_{1556}(\bv) \\
&\quad \textstyle - \frac{3}{2}Q_{1466}(\bv)Q_{5566}(\bu) + 3Q_{1266}(\bu)\wp_{66}(\bv) 
+ \frac{1}{2}Q_{2366}(\bu)\wp_{36}(\bv) \\
&\quad \textstyle - \frac{1}{10}Q_{2445}(\bu)\wp_{36}(\bv) + 4Q_{1466}(\bu)\wp_{44}(\bv) 
- \frac{3}{10}\wp_{33}(\bu)Q_{3566}(\bv) \\
&\quad \textstyle + \frac{1}{30}Q_{2556}(\bu)Q_{3566}(\bv) - \frac{4}{3}Q_{1444}(\bu)\wp_{66}(\bv) 
+ \frac{32}{3}\wp_{16}(\bu)\wp_{24}(\bv) \\
&\quad \textstyle + \frac{8}{3}\wp_{16}(\bu)Q_{4455}(\bv) + \frac{106}{15}\wp_{16}(\bu)Q_{3566}(\bv) 
+ \frac{1}{20}Q_{2566}(\bu)Q_{3556}(\bv) \\
&\quad \textstyle - 6Q_{1356}(\bv)\wp_{66}(\bu) + \frac{3}{2}Q_{1556}(\bu)Q_{5566}(\bv) \big] \lambda_4^2 
\end{align*}
\begin{align*}
&P_{18} = \textstyle \big[ \frac{3}{10}Q_{2566}(\bu)\wp_{36}(\bv) - \wp_{23}(\bu)\wp_{45}(\bv) 
- \frac{1}{5}Q_{2345}(\bv) - \frac{1}{10}Q_{3344}(\bv) \\
&\quad \textstyle + \frac{2}{5}\wp_{12}(\bu) + \frac{3}{4}\wp_{56}(\bu)Q_{2366}(\bv) 
+ \frac{1}{24}Q_{4555}(\bu)\wp_{25}(\bv) + \frac{4}{3}\wp_{34}(\bu)\wp_{25}(\bv) \\
&\quad \textstyle - \frac{1}{24}\wp_{34}(\bv)Q_{4555}(\bu)  
 + \frac{1}{3}\wp_{55}(\bu)\wp_{14}(\bv) - 2\wp_{45}(\bu)\wp_{15}(\bv) 
- \frac{1}{24}\wp_{55}(\bu)Q_{3445}(\bv) \\
&\quad \textstyle + \frac{1}{6}\wp_{55}(\bu)Q_{3366}(\bv) 
 - \frac{3}{5}Q_{3556}(\bu)\wp_{36}(\bv) \big] \lambda_2
 +  \big[ \frac{7}{20}Q_{2566}(\bu)\wp_{36}(\bv) \\
&\quad \textstyle - \frac{3}{5}Q_{3556}(\bu)\wp_{36}(\bv) 
 + 2\wp_{66}(\bu)Q_{1556}(\bv) + 2\wp_{66}(\bu)Q_{1466}(\bv) + \frac{4}{3}\wp_{44}(\bv)\wp_{16}(\bu) \\
&\quad \textstyle - \frac{2}{3}\wp_{16}(\bv)Q_{5566}(\bu)
 - \frac{1}{6}\wp_{45}(\bu)\wp_{15}(\bv) + \frac{1}{3}\wp_{23}(\bu)\wp_{45}(\bv) \\
&\quad \textstyle - \frac{1}{4}Q_{3566}(\bu)\wp_{35}(\bv) \big] \lambda_4\lambda_3
 + \big[ 4\wp_{66}(\bu)Q_{1556}(\bv) + \frac{1}{10}Q_{2566}(\bu)\wp_{36}(\bv) \\
&\quad \textstyle + 8\wp_{16}(\bv)Q_{5566}(\bu) 
+ 4\wp_{66}(\bu)Q_{1466}(\bv) - \frac{64}{3}\wp_{44}(\bv)\wp_{16}(\bu) \big] \lambda_4^3 
\end{align*}
\begin{align*}
&P_{14} = \textstyle \big[  \wp_{25}(\bu)\wp_{45}(\bv) + \wp_{22}(\bu) + \frac{1}{3}Q_{3445}(\bu) + \frac{4}{3}\wp_{14}(\bv) + \frac{5}{12}Q_{3366}(\bv) \\
&\quad \textstyle - \frac{11}{20}Q_{3556}(\bu)\wp_{56}(\bv) - 6\wp_{44}(\bu)\wp_{26}(\bv) + 3Q_{5566}(\bv)\wp_{26}(\bu) 
+ \frac{18}{5}\wp_{16}(\bu)\wp_{66}(\bv)  \\
&\quad \textstyle - \frac{1}{2}Q_{4455}(\bv)\wp_{46}(\bu) + \frac{23}{12}Q_{3566}(\bv)\wp_{55}(\bu) 
- \frac{1}{24}Q_{4555}(\bu)\wp_{45}(\bv) \\
&\quad \textstyle - \frac{27}{10}\wp_{33}(\bu)\wp_{66}(\bv) - \frac{7}{3}\wp_{24}(\bu)\wp_{55}(\bv) 
- \frac{3}{5}Q_{3456}(\bu)\wp_{66}(\bv) - 2\wp_{24}(\bv)\wp_{46}(\bu) \\
&\quad \textstyle + \frac{2}{3}Q_{4455}(\bv)\wp_{55}(\bu) - 3\wp_{35}(\bu)Q_{5566}(\bv) 
- \frac{2}{3}\wp_{34}(\bu)\wp_{45}(\bv) - \frac{11}{2}\wp_{35}(\bu)\wp_{44}(\bv) \\
&\quad \textstyle + \frac{3}{10}Q_{2556}(\bu)\wp_{66}(\bv) - \frac{3}{2}Q_{3566}(\bv)\wp_{46}(\bu) 
+ \frac{1}{10}Q_{2566}(\bv)\wp_{56}(\bu) \big] \lambda_1 \\
&\quad \textstyle + \big[ \frac{1}{24}Q_{3566}(\bv)\wp_{55}(\bu)  + \frac{1}{6}\wp_{25}(\bu)\wp_{45}(\bv) 
- \frac{1}{10}Q_{3556}(\bu)\wp_{56}(\bv) \\
&\quad \textstyle + \frac{1}{6}\wp_{34}(\bu)\wp_{45}(\bv)
- \frac{1}{20}Q_{2566}(\bv)\wp_{56}(\bu) - \frac{1}{30}Q_{3366}(\bv) - \frac{1}{15}Q_{3445}(\bu)  \\
&\quad \textstyle + \frac{8}{15}\wp_{14}(\bv) \big] \lambda_4\lambda_2
- \frac{64}{3}\wp_{16}(\bv)\wp_{66}(\bu)\lambda_4^4 
+ \frac{68}{3}\wp_{16}(\bv)\wp_{66}(\bu)\lambda_4^2\lambda_3 \\
&\quad \textstyle - \frac{32}{3}\wp_{16}(\bv)\wp_{66}(\bu)\lambda_3^2
\end{align*}
\begin{align*}
&P_{10} = \textstyle \big[ - 2\wp_{24}(\bv) - \frac{5}{3}Q_{3566}(\bv)
- \frac{5}{2}\wp_{44}(\bv)\wp_{55}(\bu) + 5\wp_{36}(\bu)\wp_{56}(\bv) \\
&\quad \textstyle - \frac{1}{24}Q_{4455}(\bv) - \frac{5}{2}\wp_{35}(\bu)\wp_{66}(\bv) 
+ \frac{5}{4}Q_{5566}(\bu)\wp_{55}(\bv) \big] \lambda_0
+ \big[ 4\wp_{44}(\bv)\wp_{46}(\bu) \\
&\quad \textstyle + \frac{8}{3}\wp_{24}(\bv) + \frac{18}{5}\wp_{36}(\bu)\wp_{56}(\bv) 
+ 6\wp_{26}(\bv)\wp_{66}(\bu) - \frac{11}{2}\wp_{35}(\bu)\wp_{66}(\bv) \\
&\quad \textstyle + 6\wp_{44}(\bv)\wp_{55}(\bu) 
- \frac{3}{2}Q_{5566}(\bu)\wp_{46}(\bv) 
+ \frac{7}{3}Q_{5566}(\bu)\wp_{55}(\bv) - \frac{5}{3}Q_{3566}(\bu) \\
&\quad \textstyle - \frac{2}{3}Q_{4455}(\bu) \big]\lambda_4\lambda_1
\end{align*}
\begin{align*} &P_{6} = \textstyle \frac{16}{3}\wp_{66}(\bu)\wp_{55}(\bv)\lambda_3\lambda_1 
+ \frac{7}{12}Q_{5566}(\bv)\lambda_4\lambda_0 - 4\wp_{46}(\bu)\wp_{66}(\bv)\lambda_4^2\lambda_1 \\
&\quad \textstyle - 4\wp_{46}(\bv)\wp_{66}(\bu)\lambda_3\lambda_1 
+ \frac{5}{2}\wp_{55}(\bu)\wp_{66}(\bv)\lambda_4\lambda_0 - 6\wp_{55}(\bv)\wp_{66}(\bu)\lambda_4^2\lambda_1 \\
&\quad \textstyle - \frac{3}{2}Q_{5566}(\bv)\lambda_3\lambda_1  + \frac{4}{3}\wp_{44}(\bv)\lambda_4\lambda_0 
- \frac{16}{3}\wp_{44}(\bv)\lambda_4^2\lambda_1 - 2Q_{5566}(\bv)\lambda_4^2\lambda_1  \\
&\quad \textstyle + 3\wp_{44}(\bu)\lambda_3\lambda_1  
\end{align*}
\begin{align*}
&P_{2} = \textstyle \frac{7}{3}\wp_{66}(\bu)\lambda_4\lambda_3\lambda_1 
- \frac{3}{5}\wp_{66}(\bu)\lambda_2\lambda_1 - \frac{16}{3}\wp_{66}(\bu)\lambda_4^3\lambda_1 
+ \frac{4}{3}\wp_{66}(\bu)\lambda_4^2\lambda_0 \\
&\quad \textstyle - \frac{2}{3}\wp_{66}(\bu)\lambda_3\lambda_0 
\end{align*}

\end{document}